\documentclass[preprint,3p,12pt]{elsarticle}

\usepackage{graphicx}
\usepackage{amsmath}
\usepackage{amsfonts}
\usepackage{enumitem} 							% Resume emurate 		
\usepackage[labelfont=bf,tableposition=top]{caption} 	% Make captions blod (Figure & Table)
\usepackage{subfig}									% Subfigure package
\usepackage{epstopdf}								% Eps to pdf, now it is possible to use both PNG and EPS images
\usepackage{floatrow}
\newfloatcommand{capbtabbox}{table}[][\FBwidth]
\newsavebox{\bigleftbox}
\floatstyle{plaintop}
\restylefloat{table}
\usepackage{color}
\usepackage{mathrsfs}  
\usepackage{multirow}
\usepackage{amsmath}
\usepackage[capitalise]{cleveref}
\usepackage{titlesec}
\usepackage{etoolbox}
\usepackage{lipsum}
\usepackage{chngpage}
\usepackage{wrapfig}

\def\B#1{\mbox{\boldmath{$#1$}}}

\newcommand{\norm}[1]{|\!| #1 |\!|}
\newcommand{\jump}[1]{\mbox{$[\![ #1 ]\!]$}}
\newcommand{\avg}[1]{\mbox{$\{\!\>\!\>\!\!\!\{ #1 \}\!\>\!\>\!\!\!\}$}}
\newcommand{\dG}{\,\text{d}\Gamma}
\newcommand{\dO}{\,\text{d}\Omega}

\newcommand*\closedset[1]{%
   \,
   \vbox{%
     \hrule height 0.5pt%                  % Line above with certain width
     \kern0.25ex%                          % Distance between line and content
     \hbox{%
       \kern -0.3em%                        % Distance between content and left side of box, negative values for lines shorter than content
       \ifmmode#1\else\ensuremath{#1}\fi%  % The content, typeset in dependence of mode
       \kern -0.3em%                        % Distance between content and left side of box, negative values for lines shorter than content
     }% end of hbox
   }% end of vbox
   \,
}

\makeatletter
\newcommand*\rel@kern[1]{\kern#1\dimexpr\macc@kerna}
\newcommand*\widebar[1]{%
  \begingroup
  \def\mathaccent##1##2{%
    \rel@kern{1.2}%
    \overline{\rel@kern{-1.2}\macc@nucleus\rel@kern{-0.2}}%
    \rel@kern{-0.2}%
  }%
  \macc@depth\@ne
  \let\math@bgroup\@empty \let\math@egroup\macc@set@skewchar
  \mathsurround\z@ \frozen@everymath{\mathgroup\macc@group\relax}%
  \macc@set@skewchar\relax
  \let\mathaccentV\macc@nested@a
  \macc@nested@a\relax111{#1}%
  \endgroup
}
\makeatother

\makeatletter
\newsavebox\myboxA
\newsavebox\myboxB
\newlength\mylenA
\newcommand*\widebarQ[2][0.6]{%
    \sbox{\myboxA}{$\m@th#2$}%
    \setbox\myboxB\null% Phantom box
    \ht\myboxB=\ht\myboxA%
    \dp\myboxB=\dp\myboxA%
    \wd\myboxB=#1\wd\myboxA% Scale phantom
    \sbox\myboxB{$\m@th\overline{\copy\myboxB}$}%  Overlined phantom
    \setlength\mylenA{\the\wd\myboxA}%   calc width diff
    \addtolength\mylenA{-\the\wd\myboxB}%
    \ifdim\wd\myboxB<\wd\myboxA%
       \rlap{\hskip 0.5\mylenA\usebox\myboxB}{\usebox\myboxA}%
    \else
        \hskip -0.5\mylenA\rlap{\usebox\myboxA}{\hskip 0.5\mylenA\usebox\myboxB}%
    \fi}
\makeatother

\begin{document}

\begin{frontmatter}

\title{\large Variationally consistent mass scaling for explicit time-integration schemes\\ of lower- and higher-order finite element methods}

\author[address1]{Stein K.F. Stoter\corref{cor1}}
\ead{K.F.S.Stoter@tue.nl}
\author[address2]{Thi-Hoa Nguyen}
\ead{Nguyen@mechanik.tu-darmstadt.de}
\author[address2]{Ren\'e R. Hiemstra}
\ead{Hiemstra@mechanik.tu-darmstadt.de}
\author[address2]{Dominik~Schillinger}
\ead{Schillinger@mechanik.tu-darmstadt.de}

\cortext[cor1]{Corresponding author}

\address[address1]{Department of Mechanical Engineering, Eindhoven University of Technology, The Netherlands}
\address[address2]{Institute of Mechanics, Computational Mechanics Group, Technical University of Darmstadt, Germany}

\begin{abstract}
In this paper, we propose a variationally consistent technique for decreasing the maximum eigenfrequencies of structural dynamics related finite element formulations. Our approach is based on adding a symmetric positive-definite term to the mass matrix that follows from the integral of the traction jump across element boundaries. The added term is weighted by a small factor, for which we derive a suitable, and simple, element-local parameter choice. For linear problems, we show % factors of increase of critical time-steps of up to 2.5 without 
that our mass-scaling method produces no adverse effects in terms of spatial accuracy and orders of convergence. We illustrate these properties in one, two and three spatial dimension, for quadrilateral elements and triangular elements, and for up to fourth order polynomials basis functions. To extend the method to non-linear problems, we introduce a linear approximation and show that a sizeable increase in critical time-step size can be achieved while only causing minor (even beneficial) influences on the dynamic response.
\end{abstract}


\end{frontmatter}

\newpage

\tableofcontents

%\linenumbers

\newpage

%=====================================================
\section{Introduction}
\label{sec:intro}
%=====================================================

Explicit time-integration techniques are suitable when the time scales of interest are inherently small and for partial differential equations that are not very stiff. Industry relevant example applications are those involving contact, fragmentation and penetration, such as explosions and crashes, or quasi-static analysis of highly nonlinear processes \cite{Abaqus633,AnsysExpDyn}.
While each time step in an explicit time-integration scheme is relatively inexpensive to compute, stability requirements limit the permissible time-step size. This so-called ``critical time step'' is related to the maximum eigenfrequency of the mass-to-stiffness generalized eigenvalue problem. The specific relation depends on the time-integration technique that is employed. A typical example for undamped structural dynamics is:
\begin{align}
    \Delta t_{crit}  = \frac{2}{\sqrt{\lambda^h_{max}}} = \frac{2}{\omega^h_{max}}\,, \label{Dtcrit}
\end{align}
which holds for the central difference method as a special case of the Newmark method~\cite[p.~493]{Hughes2000Book}.
Here, $\lambda^h_{max}$ is the maximum eigenvalue of the discrete formulation and $\omega^h_{max}$ is the corresponding maximum eigenfrequency.

\begin{figure}[!b]
    \centering
    \includegraphics[trim=0 15 10 10, clip,width=0.85\textwidth]{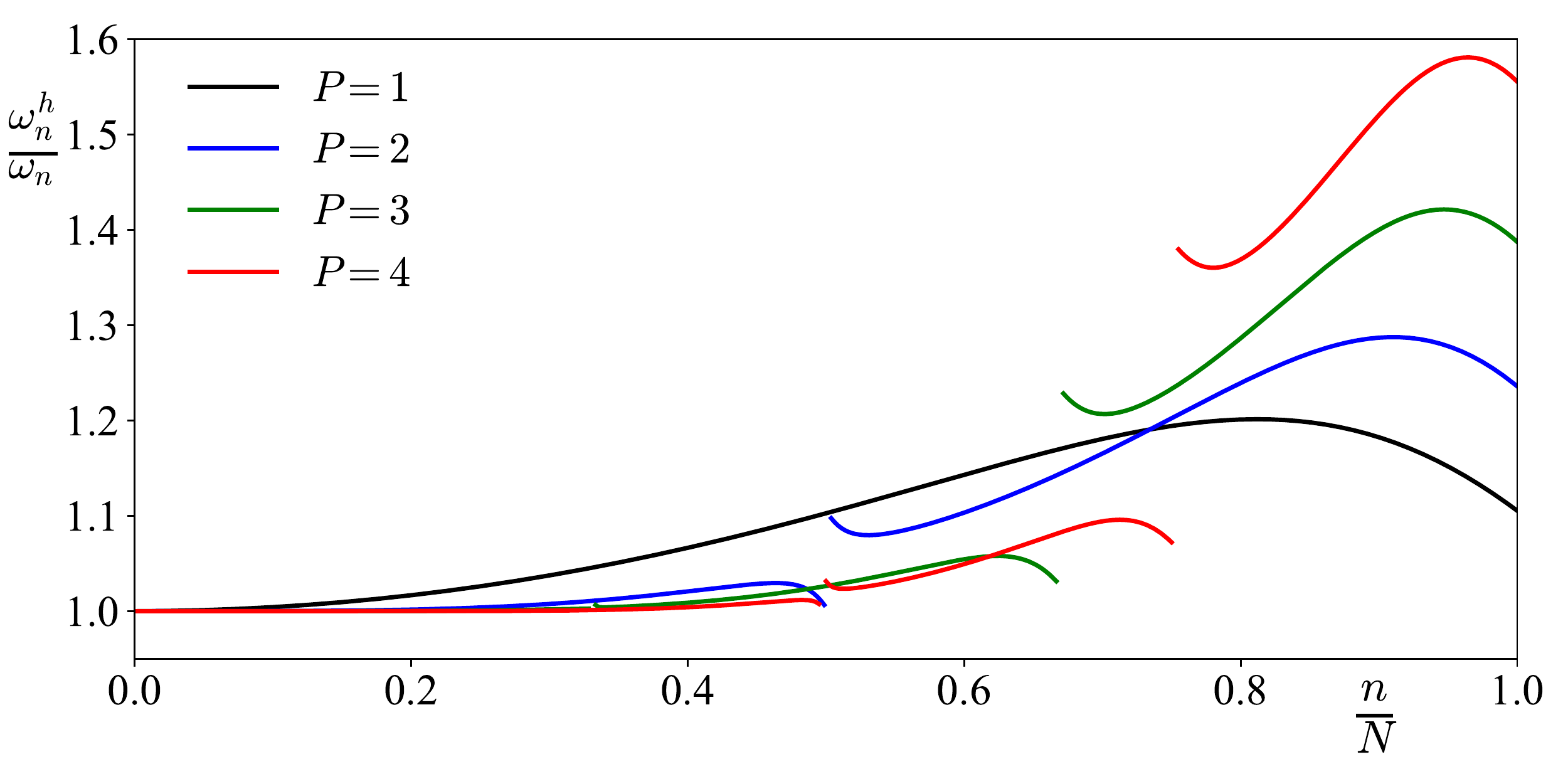}
    \caption{Relative discrete frequency spectra for a one-dimensional linear wave equation.}
    \label{fig:freqspec}
\end{figure}

\Cref{fig:freqspec} shows the well-known ratios between the discrete and analytical eigenfrequencies for a one-dimensional linear wave equation %with Dirichlet boundary conditions that is 
discretized with standard $C^0$-continuous finite elements of polynomial orders $P$ of 1 to 4. As this graph implies, the numerical eigenfrequencies of a standard Galerkin approximation always overpredict the analytical ones \cite{Hughes2000Book}. Importantly, the low-frequency part of the spectra is accurate, and they become more accurate with increasing polynomial order \cite{Hughes2000Book,Tong1971}. The spectra corresponding to the higher-order elements exhibit jumps, leading to different ``branches''. The first branch is called the acoustic branch, and is the most accurate. The other branches are called the optical branches~\cite{Brillouin1946}. The maximum frequency, which limits the critical time step, is severely overpredicted in these optical branches and diverges with increasing polynomial order.

Various techniques exist to reduce maximum eigenfrequencies. These mostly center around modifying the mass matrix. Classical examples are mass-lumping methods, which operate on the consistent mass matrix (in the sense of \cite{Archer1963,Leckie1963}) to produce a diagonal mass matrix~\cite{Hughes2000Book,Zienkiewicz:00.1}. While the primary advantage of a diagonal mass matrix is that its inversion is trivial%. Nowadays, this property is less important due to efficient LU-factorization and storage methods such that the assembly of the residual vector for complex material models typically dominates the computational cost \cite{}. However
, a secondary effect of mass lumping is a reduction of the maximal eigenvalues \cite{Tong1971}. Various lumping strategies exist, such as row-sum lumping (which may lead to negative diagonal values for higher-order serendipity elements, which are catastrophic for numerical simulation) \cite{Zienkiewicz:00.1}, diagonal scaling \cite{Hinton1976}, lumping by nodal quadrature \cite{Strang1973,Patera1984,Schillinger2014} and manifold-based methods \cite{Yang2017}. The approximations underlying these methods imply that the final formulations are not \textit{variationally consistent}, in the sense of \cite[Ch.~4]{Strang1973}. As a consequence, the true analytical solution corresponding to the strong form no-longer satisfies the discrete problem statement, resulting in sub-optimal convergence. 

Other approaches aim to add some form of additional mass to the mass matrix. Such ``mass-scaling techniques'' may or may not be variationally consistent depending on the origin of this fictitious mass. Examples are \cite{Macek1995,Olovsson2005,Olovsson2006,Askes2011}, where a weighting of some form of stiffness matrix is used as a mass scaling, \cite{Tkachuk2013}, where the addition to the mass matrix follows variationally from a penalized Hamilton's principle, and \cite{Cocchetti2013}, where nodal accelerations are modified to decrease the maximal eigenfrequency for solid shell elements.

In the current article, we explore a mass-scaling technique based on a symmetric penalization of the natural transmission conditions across element interfaces. Such techniques are effective at reducing outlier frequencies due to reduced continuity across patch interfaces and on domain boundaries in an isogeometric analysis framework \cite{Deng2021,Nguyen2022a}. Penalties on the lack of continuity are also a common theme in stabilized methods, such as ghost-penalty for immersed finite element methods \cite{Burman2010}, edge-stabilization for inf-sup stability of equal-order pairs in mixed formulations \cite{hoang2019}, and discontinuous Galerkin methods in general \cite{Cockburn2018}. In each of those cases, though, the additional terms alter the stiffness matrix rather than the mass matrix. This serves to \textit{increase} the \textit{smallest} eigenvalues. %, stabilizing the formulation and reducing the matrix condition number. 
To \textit{reduce} the \textit{largest} eigenvalues, we add a similar term to the mass matrix. %In the current work, we push this to the extreme by considering each $C^0$-finite element a single patch of an B-spline discretization \cite{}. Additionally, we rework the formulation and parameter estimation to reduce the implementational effort, ma.

The remainder of this article is structured as follows. In \cref{sec:formulation}, we propose our mass-scaling method for the linear wave equation. We also develop an estimation strategy for the involved penalty factor. In \cref{sec:spectra}, we investigate the impact of our scaled-mass formulation on discrete spectra for one- and two-dimensional domains. Afterwards, in \cref{sec:dynamics}, we perform dynamic computations and determine the impact on the convergence behavior. We then extend our formulation to a non-linear application in \cref{sec:3Dexample}, and study the impact of necessary simplifications on the solution quality. In \cref{sec:conclusion}, we draw conclusions and make suggestions for further research.

%=====================================================
\section{Derivation of the mass-scaled semi-discrete formulation}
\label{sec:formulation}
%=====================================================

We first develop our proposed finite element formulation with a variationally consistent scaled mass for the example of a linear wave equation. To get comparable results for different simulation cases, we propose a parameter estimation for the scaling factor that takes into account the material parameters, spatial dimension, mesh size and polynomial order.

\subsection{Variational form of the wave equation}
Consider a $d$-dimensional spatial domain $\Omega\subset\mathbb{R}^d$ that is partitioned into an arbitrary number of open subdomains $\Omega_i,\,\, i=1,\cdots,N$. In each subdomain, the linear wave equation reads:
\begin{align}
\rho\frac{\partial^2 }{\partial t^2}u - \nabla\cdot( T \nabla u ) &= 0 \quad \text{in } \Omega_i \times \mathcal{T}\quad i=1,\cdots,N  \,,\label{PDE}
\end{align}
where $u(t,\B{x})$ is the dependent variable, $\rho$ the positive mass per unit volume, and $T$ the positive tensile pre-stress as load per unit surface area. The temporal domain is denoted $\mathcal{T}=]0,T[$ and has end time $T$. At time $t=0$ we require the initial conditions:
\begin{subequations}\label{IC}
\begin{alignat}{5}
u &= u_{\text{IC}} \quad&& \text{on } \Omega\times\{0\} \,, \\
\frac{\partial }{\partial t} u &= v_{\text{IC}} \quad&& \text{on } \Omega\times\{0\} \,.
\end{alignat}
\end{subequations}
Additionally, we require at all times and along the entire boundary of $\Omega$, denoted $\partial\Omega$, either Dirichlet (on $\partial\Omega_D$) or Neumann (on $\partial\Omega_N$) conditions:
\begin{subequations}
\begin{alignat}{5}
u &= 0 \quad&& \text{on } \partial\Omega_D \times \mathcal{T} \,,  \label{bc1}\\
-T\nabla u \cdot \B{n} =: -T\partial_{\boldsymbol{n}} u &= g \quad&& \text{on } \partial\Omega_N \times \mathcal{T} \,,  \hspace{2.5cm} \label{bc2}
\end{alignat}
\end{subequations}
where $\B{n}$ denotes the outward facing normal vector. We assume homogeneous Dirichlet conditions purely for ease of notation later on.

Finally, the field $u(t,\B{x})$ is subjected to transmission conditions that couple the solution from subdomain to subdomain. We denote by $\Gamma$ the collection of interfaces between the subdomains $\Omega_i$. At each point on $\Gamma$, we choose a unit normal vector $\B{n}$ as either one of the outward facing unit normal vectors of the neighboring patches. Then, the appropriate transmission conditions across $\Gamma$ are:
\begin{subequations}
\begin{alignat}{5}
\jump{u} &= 0 \quad&& \text{on }\Gamma \times \mathcal{T}  \,, \label{trans1} \\
\jump{T\partial_{\boldsymbol{n}} u} &= 0 \quad&& \text{on }\Gamma  \times \mathcal{T}  \,,  \label{trans2}
\end{alignat}
\end{subequations}
where $\jump{\cdot}$ denotes the jump operator:
\begin{align}
    \jump{f} = \lim\limits_{\epsilon\rightarrow 0} \Big( f(\B{x}-\epsilon \B{n}) - f(\B{x}+\epsilon \B{n}) \Big) \,.
\end{align}
\Cref{trans1} represents the essential kinematic compatibility condition, and \cref{trans2} the natural interfacial equilibrium condition \cite{belytschko2014}.

After multiplying \cref{PDE} by a test function, integrating over $\Omega$ and performing integration by parts on each subdomain $\Omega_i$ separately, we find that a weak solution $u$ must satisfy the following statement almost everywhere in time:
\begin{align}
    \int\limits_\Omega \rho \frac{\partial^2 }{\partial t^2}u\, w \dO + \int\limits_\Omega T \nabla u \cdot \nabla w \dO - \int\limits_\Gamma \jump{T \partial_{\boldsymbol n} u} w \dG - \int\limits_{\partial\Omega}  T\partial_{\boldsymbol n} u \, w \dG = 0 \quad \forall\, w \in\mathcal{C}^\infty_0(\Omega) \,.
\end{align}
An appropriate weak formulation then follows after substitution of the natural conditions \cref{bc2,trans2}, and by searching for $u(t,\B{x})$, which, for almost every $t$, lies in the Sobolev space $H^1_0(\Omega)$, while its second time derivative is a member of the dual space $H^{-1}(\Omega)$. The space of test functions may also be widened to $H^1_0(\Omega)$, resulting in \cite[Ch.~6]{Ern2004}:
\begin{align}
&\text{Find } u \in H^1\big(\mathcal{T};H^1_0(\Omega)\big) \text{ s.t., for a.e. }t,\, \forall\,  w \in H^1_0(\Omega): \nonumber\\
&\quad \big\langle \rho\frac{\partial^2 }{\partial t^2}u , w \big\rangle_{H^{-1},H^{1}_0} + \int\limits_\Omega  T\nabla u \cdot \nabla w \dO = - \int\limits_{\partial\Omega_N} \! g \, w \dG \,,  \label{weakform}
\end{align}
where $H^1\big(\mathcal{T};H^1_0(\Omega)\big)$ is the collection of $H^1_0(\Omega)$-valued functions that satisfy the smoothness requirements of $H^1$ in the temporal domain $\mathcal{T}$, and where $\big\langle \cdot  , \cdot \big\rangle_{V',V}$ denotes a duality pairing between members of the space $V$ and its dual. 

\subsection{Finite element formulation with variationally consistent mass scaling}

To obtain a discrete formulation of the above problem, we use the method of lines. This requires a finite element approximation only in space. The usual finite element formulation of \cref{weakform} follows from Galerkin's method as:
\begin{align}
&\text{Find } u^h \in U^h \text{ s.t. } \forall\,  w^h \in U^h: \nonumber\\
&\quad \int\limits_\Omega \rho \frac{\partial^2 }{\partial t^2}u^h\, w^h \dO + \int\limits_\Omega T \nabla u^h \cdot \nabla w^h \dO  = -\int\limits_{\partial\Omega_N} \! g \, w^h \dG  \,, \label{FE1}
\end{align}
where $U^h$ is a discrete subspace of $H^1_0(\Omega)$ as any linear combination of basis functions defined on the mesh: $U^h = \text{span}\{ N_i(\B{x}) \}_{i=1}^{\!_{N_{dofs}}}$. Then, the trial and test functions can be represented as weighted sums of these basis functions, i.e.:
\begin{subequations}\label{discuh}
\begin{alignat}{2}
&u^h(t,\B{x}) &&= \sum_i^{N_{dofs}} N_i(\B{x}) \hat{u}_i^h(t) = \underbar{N}^{\text{T}}(\B{x}) \, \underbar{u}^h(t)  \,, \\
&\,  w^h(\B{x}) &&= \sum_i^{N_{dofs}} N_i(\B{x}) \hat{w}^h_i  = \underbar{N}^{\text{T}}(\B{x}) \, \underbar{w}^h \,,
\end{alignat}
\end{subequations}
and the vector of coefficients for the discrete initial conditions are based on an interpolation or projection of the fields from \cref{IC}:
\begin{subequations}\label{ICh}
\begin{alignat}{5}
&\underbar{u}^h(0)  = \mathcal{I}[ u_{\text{IC}} ]\,, \\
&\frac{\partial }{\partial t} \underbar{u}^h(0)  = \mathcal{I}[ v_{\text{IC}} ]\,.
\end{alignat}
\end{subequations}

We note that the space $U^h$ ensures that $u^h$ satisfies the \textit{essential} conditions of \cref{bc1,trans1} by construction. We now propose to add variationally consistent terms that follow from the \textit{natural} conditions of \cref{bc2,trans2} to improve the spectral properties of the formulation. We may take the second time derivative of each of these conditions to find:
\begin{subequations}
\begin{alignat}{5}
T \partial_{\boldsymbol{n}} \left(\frac{\partial^2 }{\partial t^2}u\right) &= -\ddot{g} \quad&& \text{on } \partial\Omega_N \times \mathcal{T} \,,  \label{bc2b}\\
\jump{ T \partial_{\boldsymbol{n}} \left(\frac{\partial^2 }{\partial t^2}u\right)} &= 0 \quad&& \text{on }\Gamma  \times \mathcal{T}  \,.\label{trans2b}
\end{alignat}
\end{subequations}
%where we assumed $K$ is constant in time. 

If we consider each element of the mesh a separate subdomain $\Omega_i$, then $\Gamma$ becomes the collection of interior element facets.  We then add the following two symmetrically weighted boundary integrals of the conditions in \cref{bc2b,trans2b} to the discrete formulation:
\begin{subequations}\label{bilinpenalty}
\begin{alignat}{5}
    \text{Term 1:} \qquad& \int\limits_{\partial\Omega_N}   \beta\, \dfrac{\rho}{|T|^2}  T\partial_{\boldsymbol{n}} \left(\frac{\partial^2 }{\partial t^2}u^h\right) \, T\partial_{\boldsymbol{n}} w^h \dG  + \int\limits_{\partial\Omega_N} \beta\, \dfrac{\rho}{|T|^2}\, \ddot{g} \, \partial_{\boldsymbol{n}} w^h \dG  \,, \\
    \text{Term 2:} \qquad&\int\limits_\Gamma \beta  \,\avg{\dfrac{\rho}{|T|^2}} \jump{ T\partial_{\boldsymbol{n}} \left(\frac{\partial^2 }{\partial t^2}u^h\right)}\,\jump{T\partial_{\boldsymbol{n}} w^h} \dG = 0 \,, 
\end{alignat}
\end{subequations}
with $\avg{\cdot}$ the average operator:
\begin{align}
    \avg{f} = \lim\limits_{\epsilon\rightarrow 0} \Big( \frac{1}{2}f(\B{x}-\epsilon \B{n}) + \frac{1}{2} f(\B{x}+\epsilon \B{n}) \Big) \,.
\end{align}
The additional multiplication with the physical parameters $\rho$ and $T$ in \cref{bilinpenalty} ensures dimensional consistency. In this form, $\beta$ needs to have units $m^3$, for which we introduce a mesh-size dependency later on. The second term simplifies if $\rho$ and $T$ are continuous across element edges, but the current form is also variationally consistent if the mesh conforms to an internal material interface. If $T$ is a tensor rather than a scalar (as in \cref{sec:3Dexample}) then we define $|T|$ as its largest eigenvalue.

For now, we assume a constant scalar valued $\rho$ and $T$. By adding the terms of \cref{bilinpenalty} to the finite element formulation of \cref{FE1}, we obtain the following modified form:
\begin{align}
&\text{Find } u^h \in U^h \text{ s.t. } \forall\,  w \in U^h: \nonumber\\
\begin{split}
&\int\limits_\Omega\rho\frac{\partial^2 }{\partial t^2}u^h\, w^h \dO  + \int\limits_\Gamma \beta\, \jump{ \partial_{\boldsymbol{n}} \left(\rho \frac{\partial^2 }{\partial t^2}u^h\right)}\,\jump{\partial_{\boldsymbol{n}} w^h} \dG + \int\limits_{\partial\Omega_N}   \beta\,  \partial_{\boldsymbol{n}} \left(\rho \frac{\partial^2 }{\partial t^2}u^h\right) \, \partial_{\boldsymbol{n}} w^h \dG \\
& \hspace{2cm} + \int\limits_\Omega T \nabla u^h \cdot \nabla w^h \dO   = - \int\limits_{\partial\Omega_N} \! g \, w^h \dG - \int\limits_{\partial\Omega_N} \beta\, \rho T^{-2}\, \ddot{g} \, \partial_{\boldsymbol{n}} w^h \dG \,.  \label{FEM}
\end{split}
\end{align}
%where $\beta$ is a new parameter that needs to be chosen.

After substituting the representation of the test and trial functions per \cref{discuh} into \cref{FEM}, we can rewrite it as a matrix system of equations:
\begin{align}
    ( \mathbf{M} + \beta \mathbf{M}_\Gamma )\, \ddot{\underbar{u}}^h + \mathbf{K}\, \underbar{u}^h = \underbar{F}  \,. \label{FEmat}
\end{align}
As this equation illustrates, the additional terms in the weak formulation only affect the mass matrix. The added matrix $\mathbf{M}_\Gamma$ is symmetric positive semidefinite and serves to suppress the highest eigenvalues \cite{Nguyen2022a,li_perturbation_2014}. The variational consistency of the added terms ensures that the true solution also satisfies \cref{FEM}, which is an important property for ensuring optimal convergence \cite{Strang1973,Ern2004}.

\subsection{Parameter estimation}

% The time-step limitation for explicit finite difference treatment of the time-derivative relates to the maximum frequency of the following generalized eigenvalue problem:
% \begin{align}
%     \mathbf{K}\, \underbar{u}^h = \lambda\, ( \mathbf{M} + \beta \, \mathbf{M}_\Gamma )\, \underbar{u}^h
% \end{align}
The dynamic response of the solution advanced in time by \cref{FEmat} can be studied by examining the corresponding generalized eigenvalue problem. Of particular importance is the \textit{maximum} eigenvalue, as this limits the critical time-step size for explicit time-integration algorithms. The generalized eigenvalue problems for the maximum eigenvalues with and without the mass-scaling term read:
\begin{subequations}\label{eigenvalueprob}
\begin{alignat}{2}
    &  \mathbf{K}\, \bar{\underbar{$\xi$}}  = \bar{\lambda}_{\text{max}} \mathbf{M} \bar{\underbar{$\xi$}} \label{eigv_unpert}  \,, \\
    & \mathbf{K}\, \underbar{$\xi$}  = \lambda^*_{\text{max}} ( \mathbf{M} + \beta \mathbf{M}_\Gamma )\underbar{$\xi$}  \,, \label{eigv_pert}
\end{alignat}
\end{subequations}
where $\lambda^*_{\text{max}}$ is the target value that we aim to obtain by properly choosing $\beta$. We write this target as a factor of the true eigenvalue: $\lambda^*_{\text{max}} = a\lambda_{\text{max}}$. Similarly, we write the eigenvalue of the system without mass scaling as $\bar{\lambda}_{\text{max}} = b\lambda_{\text{max}}$. Premultiplying \cref{eigv_pert} by the eigenvector from \cref{eigv_unpert} and making use of the symmetry of the matrices then gives:
\begin{align}
\bar{\underbar{$\xi$}}^\text{T} \mathbf{K}\, \underbar{$\xi$}  = a\lambda_{\text{max}} \bar{\underbar{$\xi$}}^\text{T} ( \mathbf{M} + \beta \mathbf{M}_\Gamma )\underbar{$\xi$} = \frac{a}{b} \bar{\underbar{$\xi$}}^\text{T} \mathbf{K}\underbar{$\xi$}  + a\lambda_{\text{max}} \beta \bar{\underbar{$\xi$}}^\text{T} \mathbf{M}_\Gamma \underbar{$\xi$}  \,, 
\end{align}
or:
\begin{align}
\beta = \left(\frac{1}{a}-\frac{1}{b}\right) \frac{1}{\lambda_{\text{max}}} \,\frac{  \bar{\underbar{$\xi$}}^\text{T} \mathbf{K}\underbar{$\xi$}   }{ \bar{\underbar{$\xi$}}^\text{T} \mathbf{M}_\Gamma \underbar{$\xi$}  } \,. \label{lmbdbefore}
\end{align}
To obtain a closed and element-local estimation strategy for $\beta$ we:
\begin{enumerate}
    \item Choose an $a$ to eliminate the unknown $b$ from $\left(\frac{1}{a}-\frac{1}{b}\right)$.
    \item Approximate $\lambda_{\text{max}}$ in terms of local mesh size $h$ and polynomial order $P$.\\[-0.75cm]
    \item Estimate $\dfrac{  \bar{\underbar{$\xi$}}^\text{T} \mathbf{K}\,\underbar{$\xi$}   }{ \bar{\underbar{$\xi$}}^\text{T} \mathbf{M}_\Gamma \, \underbar{$\xi$}  }$ based on a simplified problem.
\end{enumerate}

\subsubsection{Choice of critical time step increase}
\label{ssec:abc}
First, we choose $a$, which represents the fraction between the maximum eigenvalue from \cref{eigv_pert} and the true eigenvalue. We base our choice on a target improvement of the critical time step. The critical time step scales inversely to the square root of the eigenvalue. So, for increasing the critical time-step size by a factor $\vartheta$, we require:
\begin{align}
    \sqrt{\frac{b}{a}} = \vartheta \,. \label{abC}
\end{align}
We choose the following expression for $\vartheta$:
\begin{align}
    \vartheta = \sqrt{cb+1} \,. \label{C}
\end{align}
For a fixed $c>0$, the factor $\vartheta$ is always larger than 1 and increases with $b$. This means that the critical time step will always increase due to the scaled mass, and it will increase more significantly when the discrete formulation without mass scaling severely overpredicts the largest eigenvalues. 
For a larger $c$, the time step is increased more, but this may come at the cost of decreased accuracy of the higher frequency response. 

The choice of \cref{abC} in combination with \cref{C} also simplifies the expression for $\beta$:
\begin{align}
     \left( \frac{1}{a} - \frac{1}{b}\right) = \left( \frac{cb+1}{b} - \frac{1}{b}\right) = c \label{rerabc} \quad \text{so that} \quad \beta = c \frac{1}{\lambda_{\text{max}}} \,\frac{  \bar{\underbar{$\xi$}}^\text{T} \mathbf{K}\underbar{$\xi$}   }{ \bar{\underbar{$\xi$}}^\text{T} \mathbf{M}_\Gamma \underbar{$\xi$}  }\,.
\end{align}
Based on \cref{C} and the $b$-values observed in \cref{fig:freqspec} (which range between 1.1 and 1.6), a value $c = 1$ should roughly corresponding to a $50\%$ increase in critical time step, and $c = 5$ to a 150\%-200\% increase. % We consider this the base-line value and study the effect of increase increases and decrea
%Throughout this article, we will aim for (equalling a cost reduction of $33\%$ in terms of the total number of time-steps). 

\subsubsection{Element local maximum eigenvalue estimate}
\label{ssec:trueeig}
Next, we approximate $\lambda_{\text{max}}$ by considering the analytical values on simplified domains. On a one-dimensional domain of length $L$, a two-dimensional square domain with sides $L$ and a three-dimensional cubic domain with sides $L$ the true eigenvalues are:
\begin{subequations}\label{eigenvalues}
\begin{alignat}{2}
    & \text{In 1D:}&&\quad \lambda_n = \frac{T}{\rho}\frac{\pi^2}{L^2} n^2  \,, \label{eigenvalues1D}\\
    & \text{In 2D:}&&\quad \lambda_{m,n} = \frac{T}{\rho} \frac{\pi^2}{L^2} (m^2+n^2) \,, \label{eigenvalues2D}\\
    & \text{In 3D:}&&\quad \lambda_{l,m,n} = \frac{T}{\rho} \frac{\pi^2}{L^2} (l^2+m^2+n^2) \,.
\end{alignat}
\end{subequations}
When these domains are discretized in the $x$, $y$ and $z$ directions with $n$ equidistant elements of polynomial order $P$, then the maximum eigenvalue is the $(Pn)^d$-th eigenvalue. For that maximum eigenvalue, \cref{eigenvalues} may be written in terms of the element diameter $h$:
\begin{subequations}\label{eigenvalues_num}
\begin{alignat}{2}
    & \text{In 1D:}&&\quad  \lambda_{\text{max}} = \frac{T}{\rho}\frac{\pi^2}{L^2} (Pn)^2 = \pi^2 \frac{T}{\rho} \frac{P^2}{h_x^2} = \pi^2 \frac{T}{\rho} \frac{P^2}{h^2} \,, \\
    & \text{In 2D:}&&\quad  \lambda_{\text{max}} = \frac{T}{\rho}\frac{\pi^2}{L^2} ((Pn)^2+(Pn)^2) = 2 \pi^2 \frac{P^2}{h_x^2} = 4\pi^2 \frac{T}{\rho} \frac{P^2}{h^2} \,, \\
    & \text{In 3D:}&&\quad  \lambda_{\text{max}} = \frac{T}{\rho}\frac{\pi^2}{L^2} ((Pn)^2+(Pn)^2+(Pn)^2) = 3 \pi^2 \frac{T}{\rho} \frac{P^2}{h_x^2} = 9\pi^2 \frac{T}{\rho} \frac{P^2}{h^2}  \,, 
\end{alignat}
\end{subequations}
and thus in general:
\begin{align}
    \lambda_{\text{max}} = \frac{T}{\rho} \left( d \, \pi\,  \frac{P}{h} \right)^2 \,. \label{reslmbdmax}
\end{align}

\subsubsection{Estimate of stiffness to scaled mass ratio}
Finally, we estimate the ratio between $\bar{\underbar{$\xi$}}^\text{T} \mathbf{K}\underbar{$\xi$}$ and  $\bar{\underbar{$\xi$}}^\text{T} \mathbf{M}_\Gamma \underbar{$\xi$}  $. The largest eigenvectors of \cref{eigenvalueprob} conflict most with the continuity requirement of \cref{trans2}, and thereby correspond to the largest $\bar{\underbar{$\xi$}}^\text{T} \mathbf{M}_\Gamma \underbar{$\xi$}$ \cite{Horger2019,Nguyen2022a}. Based on this empirical observation, we propose the following approximation:
\begin{align} \label{maxfunc}
\frac{  \bar{\underbar{$\xi$}}^\text{T} \mathbf{K}\underbar{$\xi$}   }{ \bar{\underbar{$\xi$}}^\text{T} \mathbf{M}_\Gamma \underbar{$\xi$}  } \approx \min\limits_{\underbar{v}\neq\underbar{0}} \frac{  \underbar{v}^\text{T} \mathbf{K}\underbar{v}   }{ \underbar{v}^\text{T} \mathbf{M}_\Gamma \underbar{v}  } = \min\limits_{0\neq v \in U^h} \frac{ \int_\Omega T \nabla v \cdot \nabla v   }{ \int_\Gamma \rho \jump{\nabla v}\jump{\nabla v} }  \,, 
\end{align} 
which we can solve analytically for a simplified case.

\begin{figure}[!t]\vspace{-0.75cm}
    \centering
    \subfloat[$P=1$]{\includegraphics[width=0.46\linewidth]{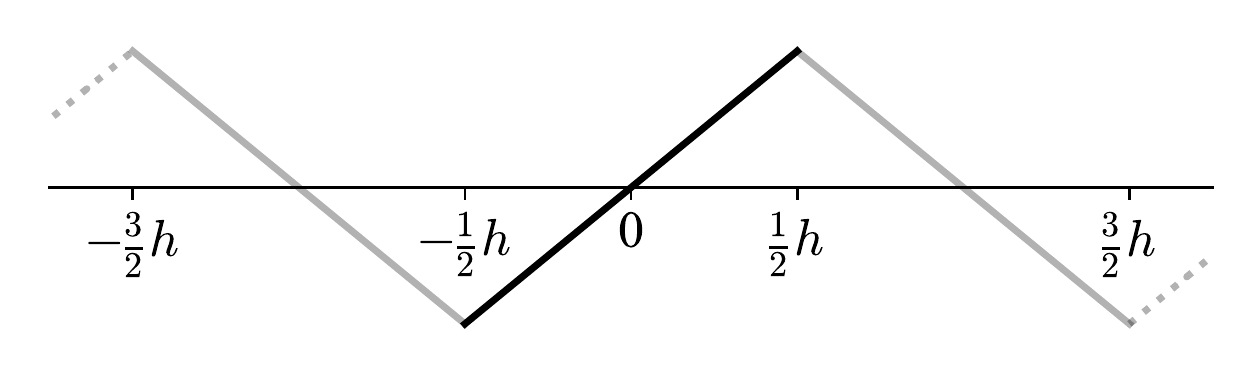}}\hspace{.25cm}
    \subfloat[$P=2$]{\includegraphics[width=0.46\linewidth]{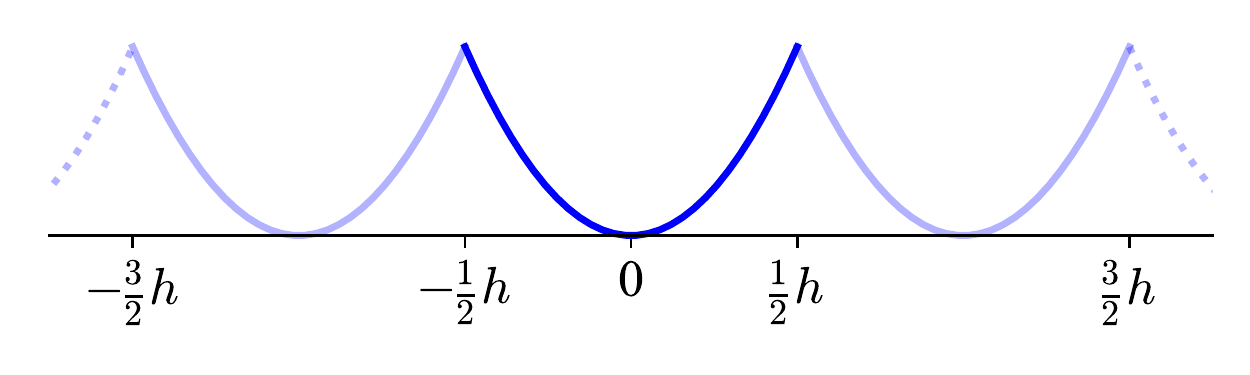}}\\
    \subfloat[$P=3$]{\includegraphics[width=0.46\linewidth]{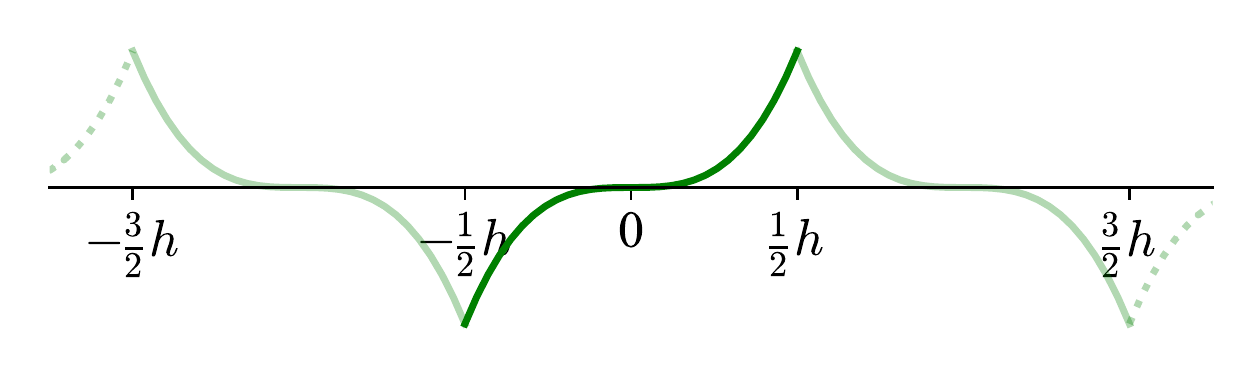}}\hspace{.25cm}
    \subfloat[$P=4$]{\includegraphics[width=0.46\linewidth]{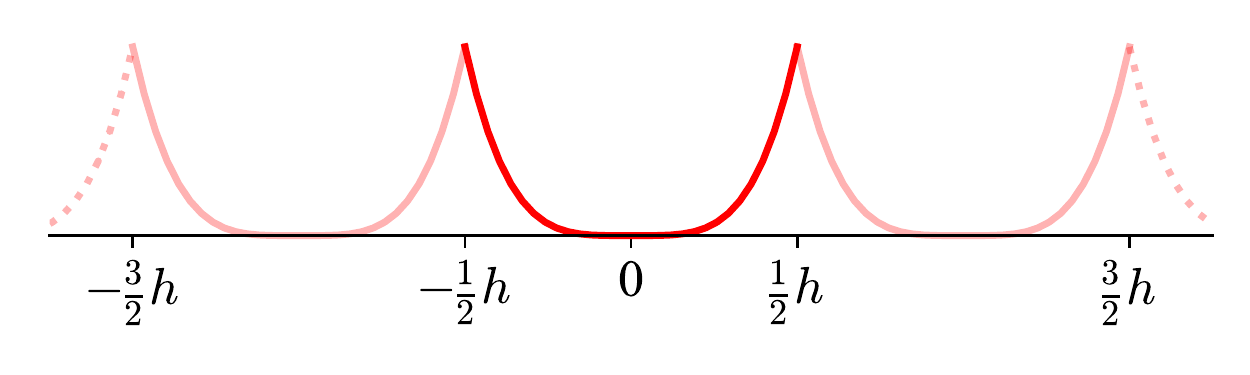}}
    \caption{Functions on a one-dimensional periodic mesh that satisfy the minimum of \cref{maxfunc}.} \label{fig:minmaxfuncs}
\end{figure}

On a one-dimensional periodic domain, the functions that maximize $\jump{\nabla v}^2$ across the element edges while minimizing $\norm{\nabla v}^2$ are depicted in \cref{fig:minmaxfuncs}. On a reference domain $[-\frac{1}{2}h,\frac{1}{2}h]$ they simply read:
\begin{align}
    v = \pm x^P \,.
\end{align}
For such a function, we can work out both the numerator and denominator as:
\begin{subequations}
\begin{alignat}{2}
 & \int_\Omega T \nabla v \cdot \nabla v \dO &&=  T\, N_{els} \int\limits_{-\frac{1}{2}h}^{\frac{1}{2}h} (P\,x^{P-1})^2 = T\, N_{els} \, \frac{2 P^2}{2P-1}\Big(\frac{h}{2}\Big)^{2P-1} \,, \\ %N_{els} \frac{P^2}{2P-1}\Big[ x^{2P-1} \Big]^{\frac{1}{2}h}_{-\frac{1}{2}h} =
 & \int_\Gamma \rho \jump{ \nabla v}\jump{ \nabla v} \dG &&= \rho \, N_{els} \, ( 2 P x^{P-1} )^2 \Big|_{\frac{1}{2}h} = \rho \, N_{els} \, 4 P^2 (\frac{h}{2})^{2P-2}  \,, 
 \end{alignat}
\end{subequations}
such that we obtain the approximation:
\begin{align}
    \frac{  \bar{\underbar{$\xi$}}^\text{T} \mathbf{K}\underbar{$\xi$}   }{ \bar{\underbar{$\xi$}}^\text{T} \mathbf{M}_\Gamma \underbar{$\xi$}  } \approx  \frac{T}{\rho}\frac{1}{4}\frac{1}{2P-1} h  \,. \label{resfrac}
\end{align}

Collecting \cref{lmbdbefore,reslmbdmax,rerabc,resfrac} finally leads to the following estimation for $\beta$:
\begin{align}
\beta = c \frac{1}{4}\frac{1}{d^2 \pi^2 }\frac{1}{2P^3-P^2} h^3 \label{beta} \,.
\end{align}

%=====================================================
\section{Effect on frequency spectra}
\label{sec:spectra}
%=====================================================

Before we study the impact of the mass-scaling term on computations of transient response, we first focus on its effect on the discrete frequency spectra. We consider a one-dimensional case and a two-dimensional case with a square domain, for both of which the true eigenvalues are known.

\subsection{The one-dimensional case: a string}
\label{ssec:1D}

\begin{figure}[!b]
    \centering
    \subfloat[$P=1$]{\includegraphics[width=0.5\linewidth]{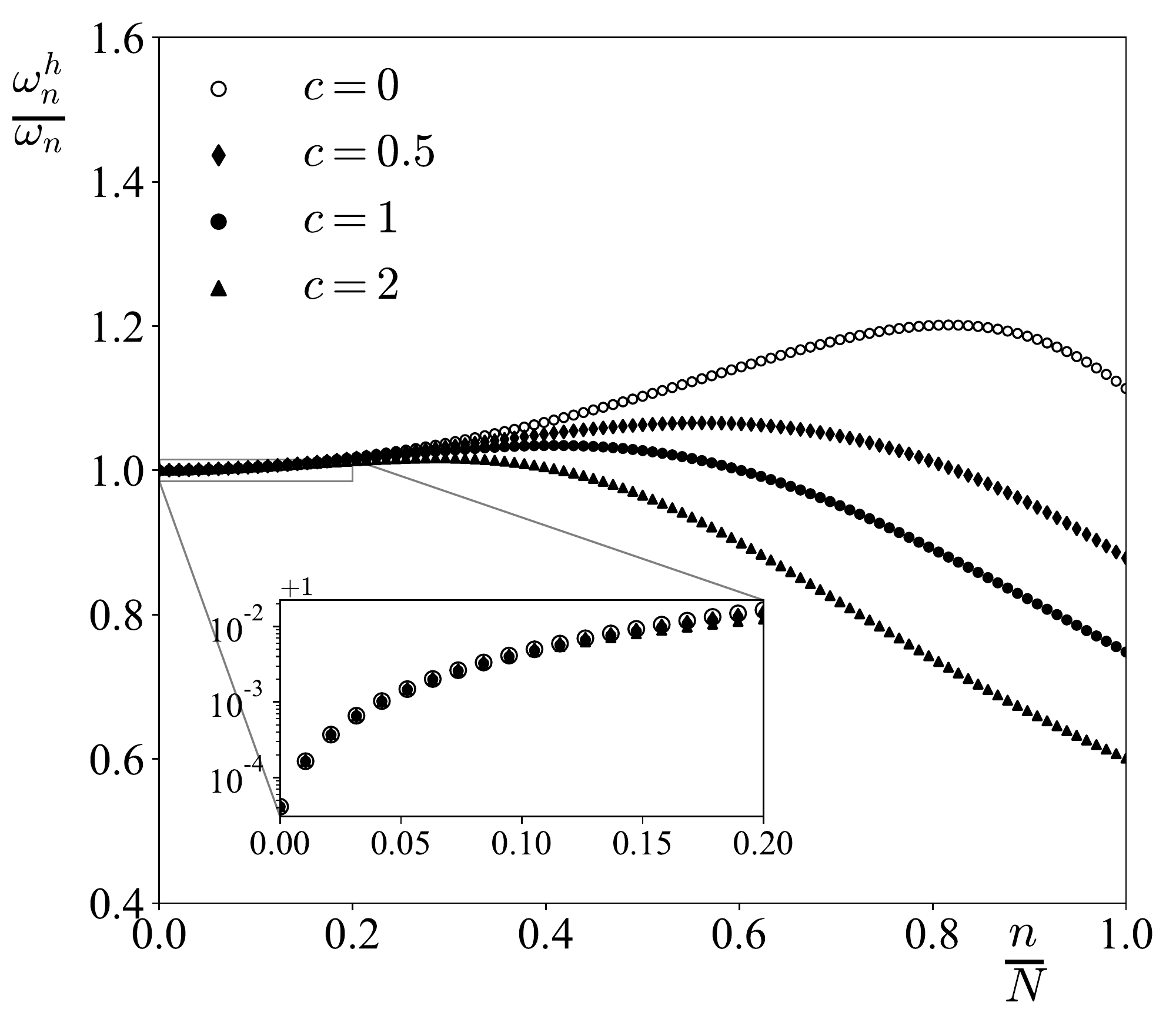}\label{fig:1D_rel1}}\hfill
    \subfloat[$P=2$]{\includegraphics[width=0.5\linewidth]{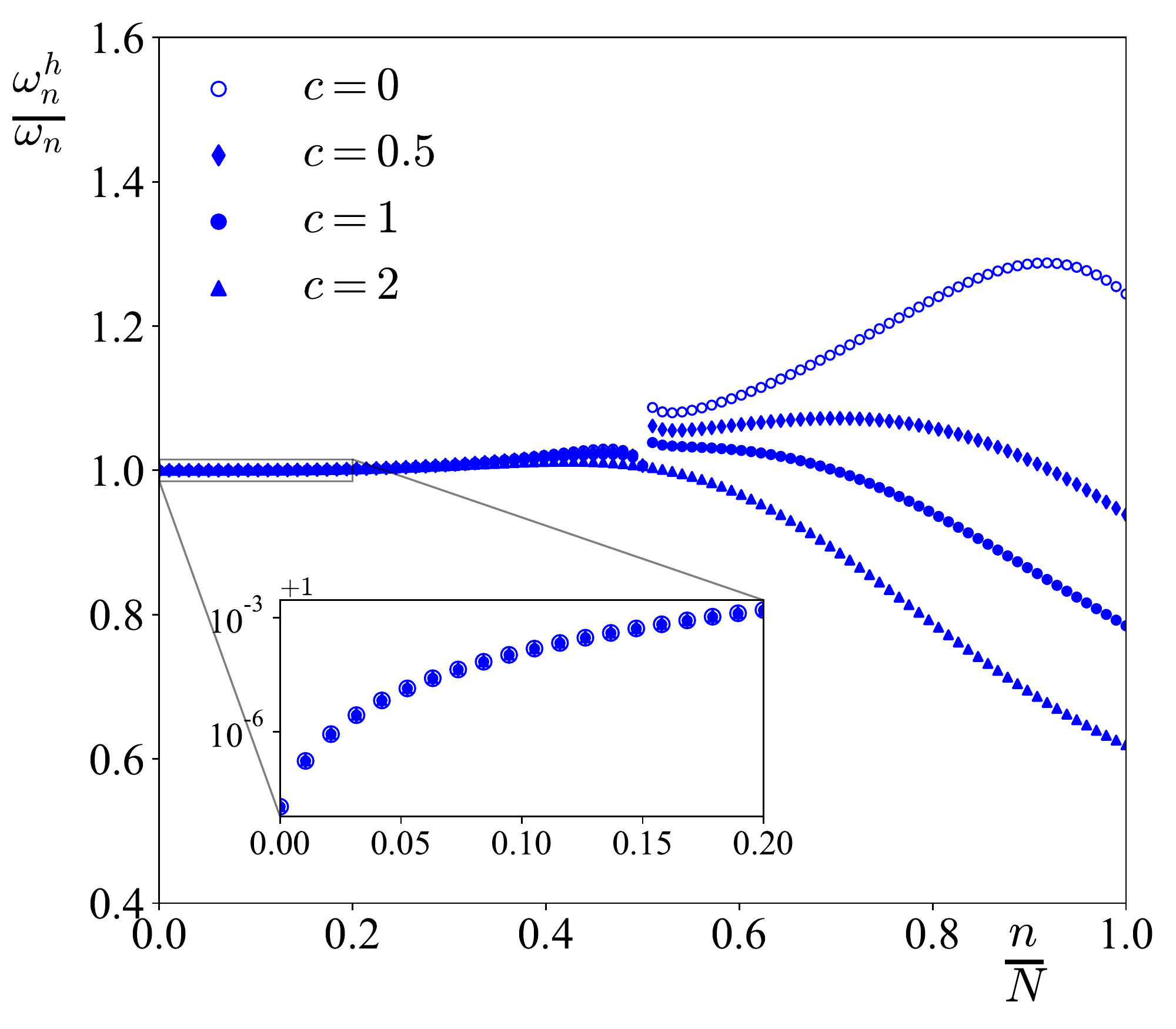}\label{fig:1D_rel2}}\\
    \subfloat[$P=3$]{\includegraphics[width=0.5\linewidth]{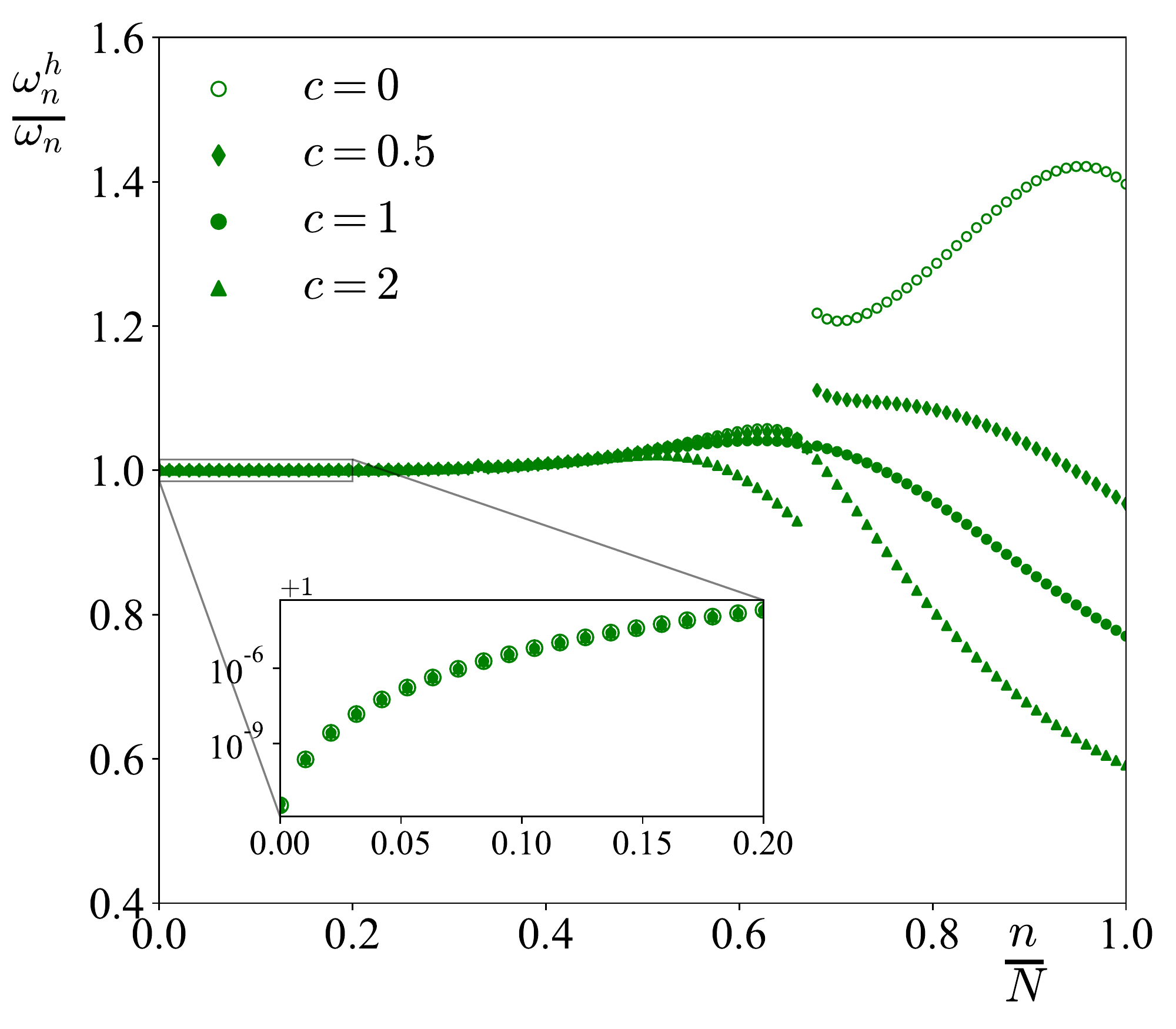}\label{fig:1D_rel3}}\hfill
    \subfloat[$P=4$]{\includegraphics[width=0.5\linewidth]{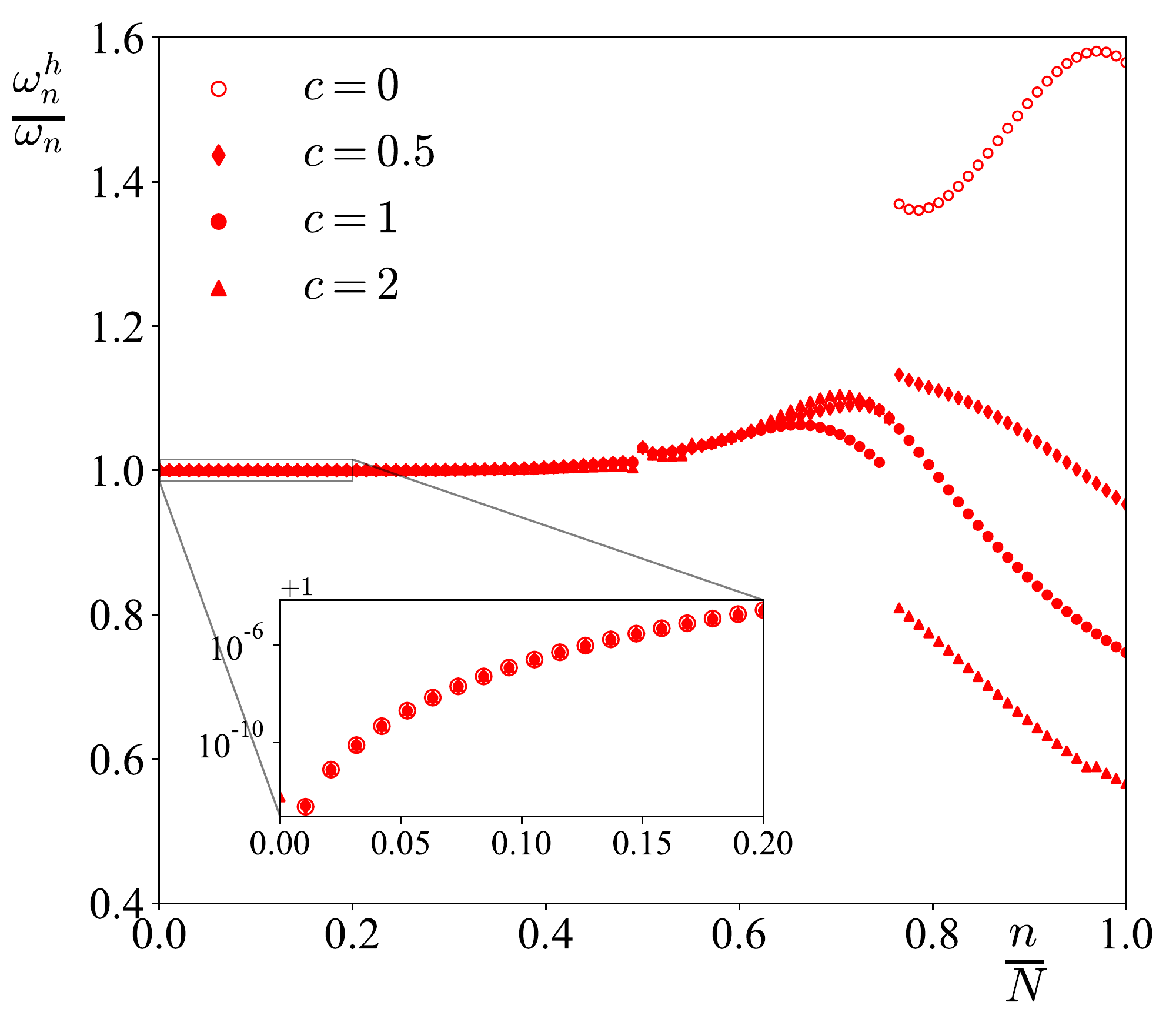}\label{fig:1D_rel4}}
    \caption{Frequency spectra for a one-dimensional domain and various polynomial orders, normalized by the exact eigenfrequencies.} \label{fig:1D_rel}
\end{figure}

We begin by analyzing the discrete frequency spectrum for the one-dimensional case with Dirichlet boundary conditions on both ends. The eigenfrequencies are the square roots of the eigenvalues:
\begin{align}
    \omega^h_n = \sqrt{\lambda^h_n} \,,
\end{align}
and the true eigenvalues are those of \cref{eigenvalues1D}. We plot the ratio between the numerically computed eigenfrequencies and the analytical eigenfrequencies as a function of the normalized mode number for polynomial orders $P=1$ until $P=4$ in \cref{fig:1D_rel1,fig:1D_rel2,fig:1D_rel3,fig:1D_rel4}. For all computations, the total number of degrees of freedom, $N$, is kept at approximately 100. The open markers correspond to the formulation without mass scaling (i.e., $c=0$) and are those of \cref{fig:freqspec}. The solid markers correspond to the eigenvalue problem of \cref{FEM,eigv_pert} which includes mass scaling for different values of $c$. The ordering of the numerical eigenfrequencies is determined based on an iterative scheme: starting with $n=0$, we select the numerical eigenfrequency that corresponds to the numerical eigenfunction that yields the lowest $L^2$-error compared to the $n$-th analytical eigenfunction. This numerical eigenfunction is removed from the selection set, and the process is repeated for the subsequent $n$. 

We observe a significant improvement of the spectrum for all polynomial orders. For all choices of $c$, the eigenfrequencies in the high frequencies range, the last $25\%$ of the spectrum, are suppressed. In many cases, the result is an underprediction of the analytical values, which is not a problem: the high frequencies do generally not contribute to the fidelity of dynamics simulations \cite[Ch.~7 and 8]{Hughes2000Book}, but do limit their critical time step. In the lowest $25\%-50\%$ of the spectrum we observe that the mass scaling does not negatively affect the numerically computed eigenfrequencies, as can be seen clearly in the semi-log inset figures. This is crucial, since these eigenfunction-eigenfrequency pairs typically dominate the system response. For the intermediate range of frequencies, i.e. until roughly $75\%$, the normalized numerical eigenfrequencies are pushed closer towards the desired value of 1, improving the overall accuracy and smoothness of the discrete spectra.

\begin{figure}[!b]
    \centering
    \subfloat[$P=1$]{\includegraphics[width=0.5\linewidth]{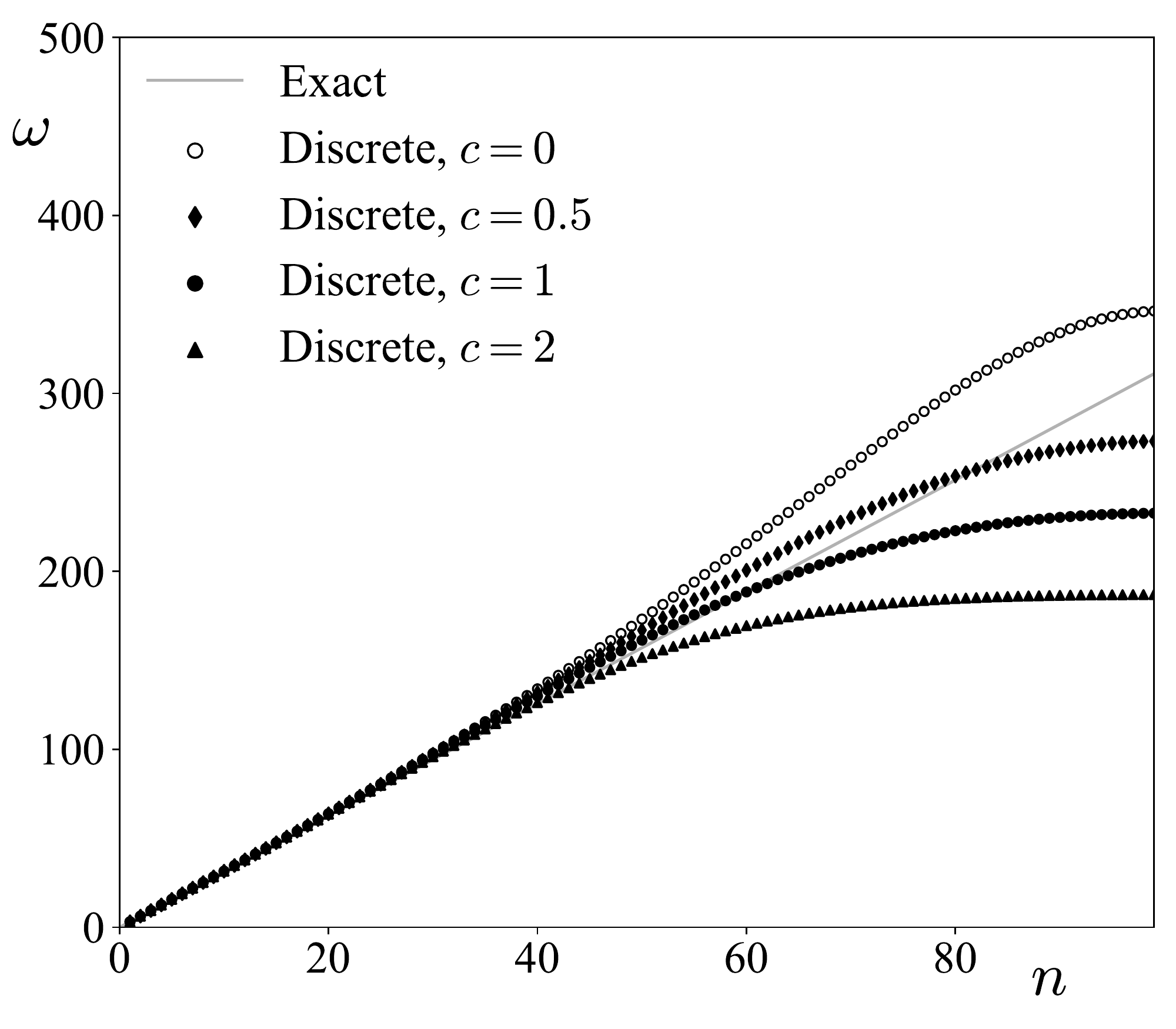}}\hfill
    \subfloat[$P=2$]{\includegraphics[width=0.5\linewidth]{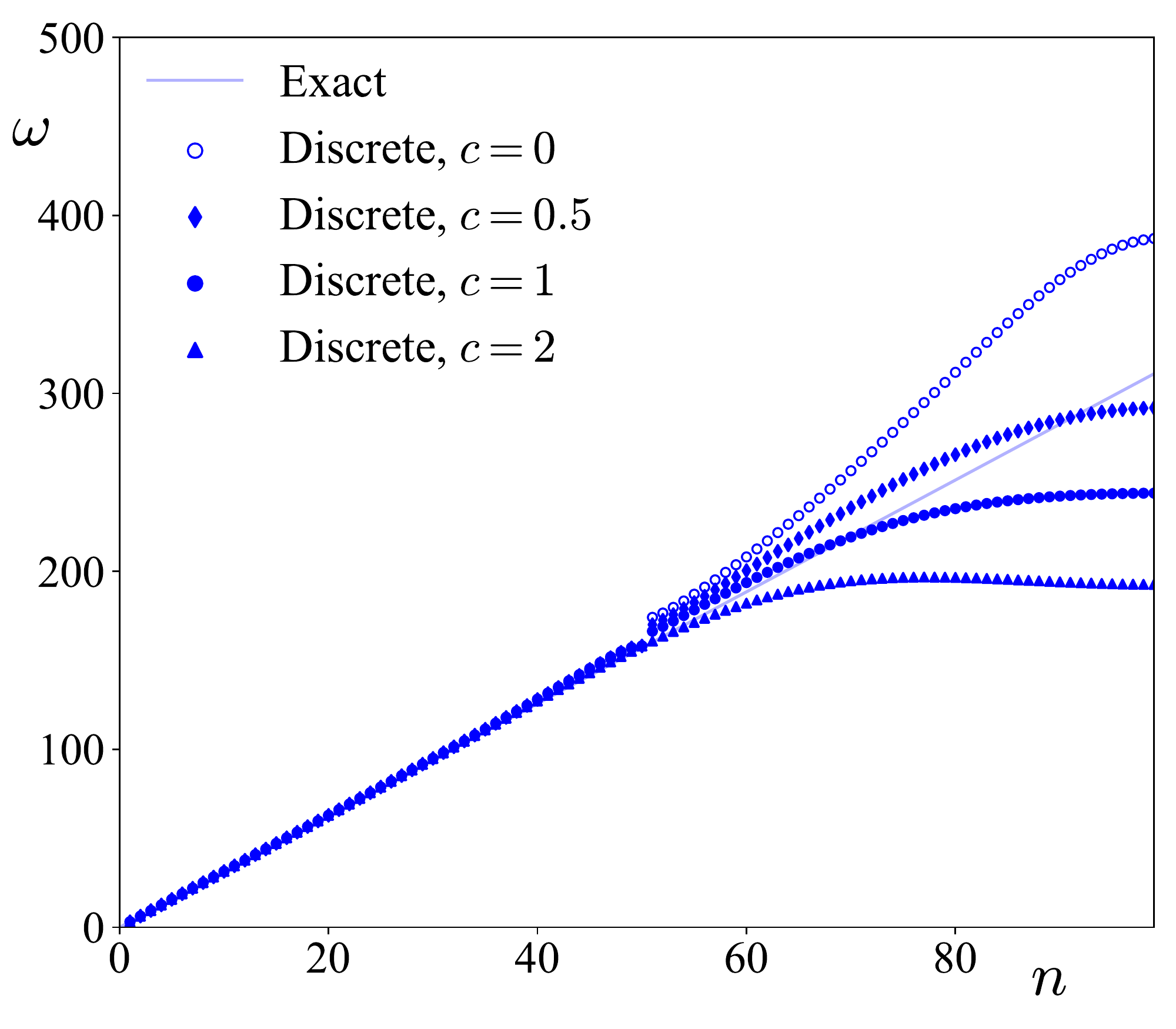}}\\
    \subfloat[$P=3$]{\includegraphics[width=0.5\linewidth]{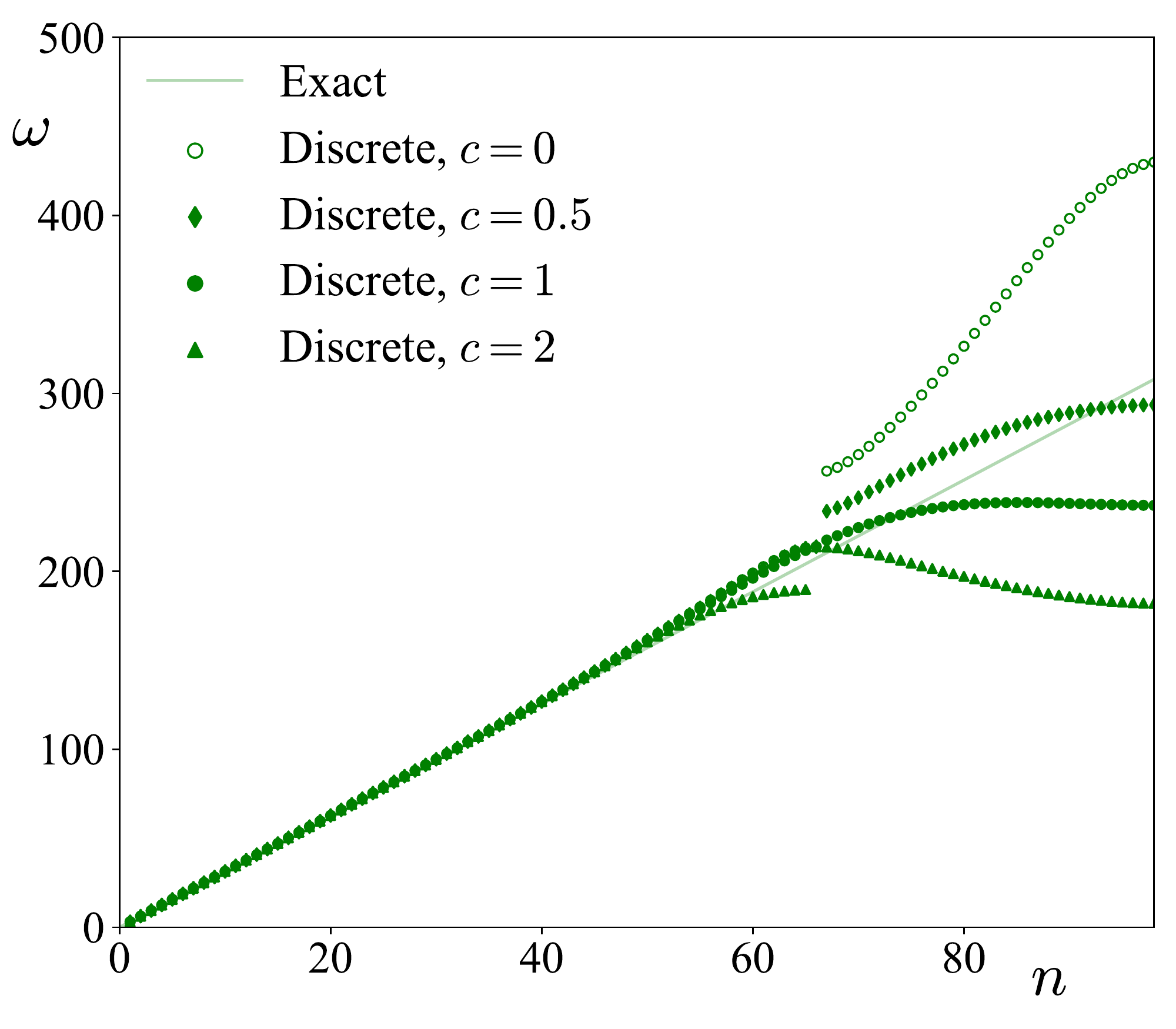}}\hfill
    \subfloat[$P=4$]{\includegraphics[width=0.5\linewidth]{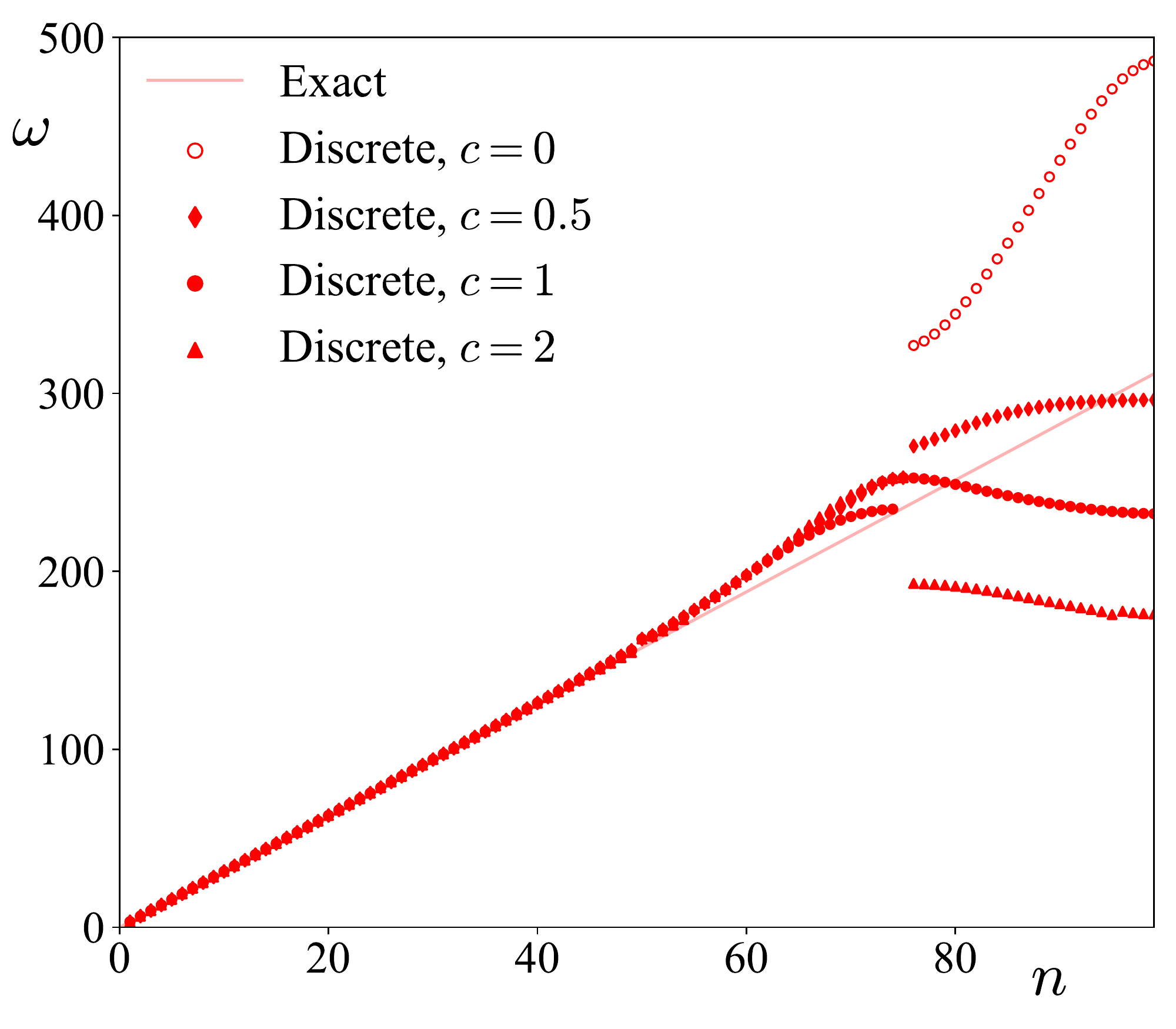}}
    \caption{Eigenfrequencies for a one-dimensional domain and various polynomial orders.} \label{fig:1D_abs}
\end{figure}

%, which is the desired result as the formulation without scaled mass already accurately predicts the frequencies in this range. In the high frequency range, the last $25\%$, the stabilized formulation under-predicts the frequencies and produces values significantly smaller than those of the non-stabilized formulation. This 

Since \cref{fig:1D_rel} shows the ratio between the numerical and analytical eigenfrequencies, the results may give the false impression that an increase in $c$ steadily decreases the \textit{maximum} eigenfrequency. In \cref{fig:1D_abs} we show the non-normalized eigenfrequencies. For $P=1$, we do observe that the maximum eigenfrequency is decreased at a steady rate for increasing $c$, but the same is not true for the cases of $P=2$, $3$ and $4$. Past a threshold value for $c$ the high frequencies are suppressed below the values of an intermediate frequency such that this become the new maximum. The intermediate frequencies only slowly decrease as $c$ becomes larger. This is particularly clear for $P=4$, where the eigenfrequencies in the before-last optical branch become dominant.

For higher-order $C^0$ finite elements, the lack of inter-element continuity is the cause of the different optical branches, and the final branch is due to the lack of first-order continuity (i.e., $U^h\subset C^0 \setminus C^1$) \cite{Brillouin1946}. This is the reason for the superior spectral properties of the higher-order continuous spline basis functions used in isogeometric analysis \cite{Cottrell:09.1}. The mass-scaling term of \cref{bilinpenalty} acts on this first-order inter-element continuity. It thus specifically targets the highest frequency branch. One potential approach for suppressing the frequencies from the lower optical branches is thus to add additional mass-scaling terms that involve jumps of higher-order derivatives.  This concept is explored in an isogeometric analysis framework across patch-interfaces in \cite{Nguyen2022a}, and on domain boundaries in \cite{Deng2021}. However, in the current article we aim to maintain the variational consistency of the finite element formulation, which only permits us to add a penalty on the natural conditions of \cref{bc2,trans2}, i.e., the lowest derivative.

%The three different curves for the stabilized formulation indicate that the results are rather sensitive to the choice of $\beta$, represented now in terms of $c$. The proposed medium $c=1$ appears to produce the optimal results in term of spectrum smoothness and decrease of the higher frequencies. The bump a the end of the spectrum for $P=4$ and $c=2$ is due to a re-ordering of the modes: with this $c$ the frequencies in the third optical branch are suppressed so far that those at the end of the second optical branch become the largest frequencies and shift to the right. The consistency in the results for various polynomial orders is remarkable, given the sensitivity to $\beta$ and the predicted non-trivial dependency on $P$ from \cref{beta}. 

\subsection{A two-dimensional case: a square drum}

\begin{figure}[!b]\vspace{-0.25cm}
    \centering
    \subfloat[$P=1$]{\includegraphics[width=0.5\linewidth,trim=0 13 0 10, clip]{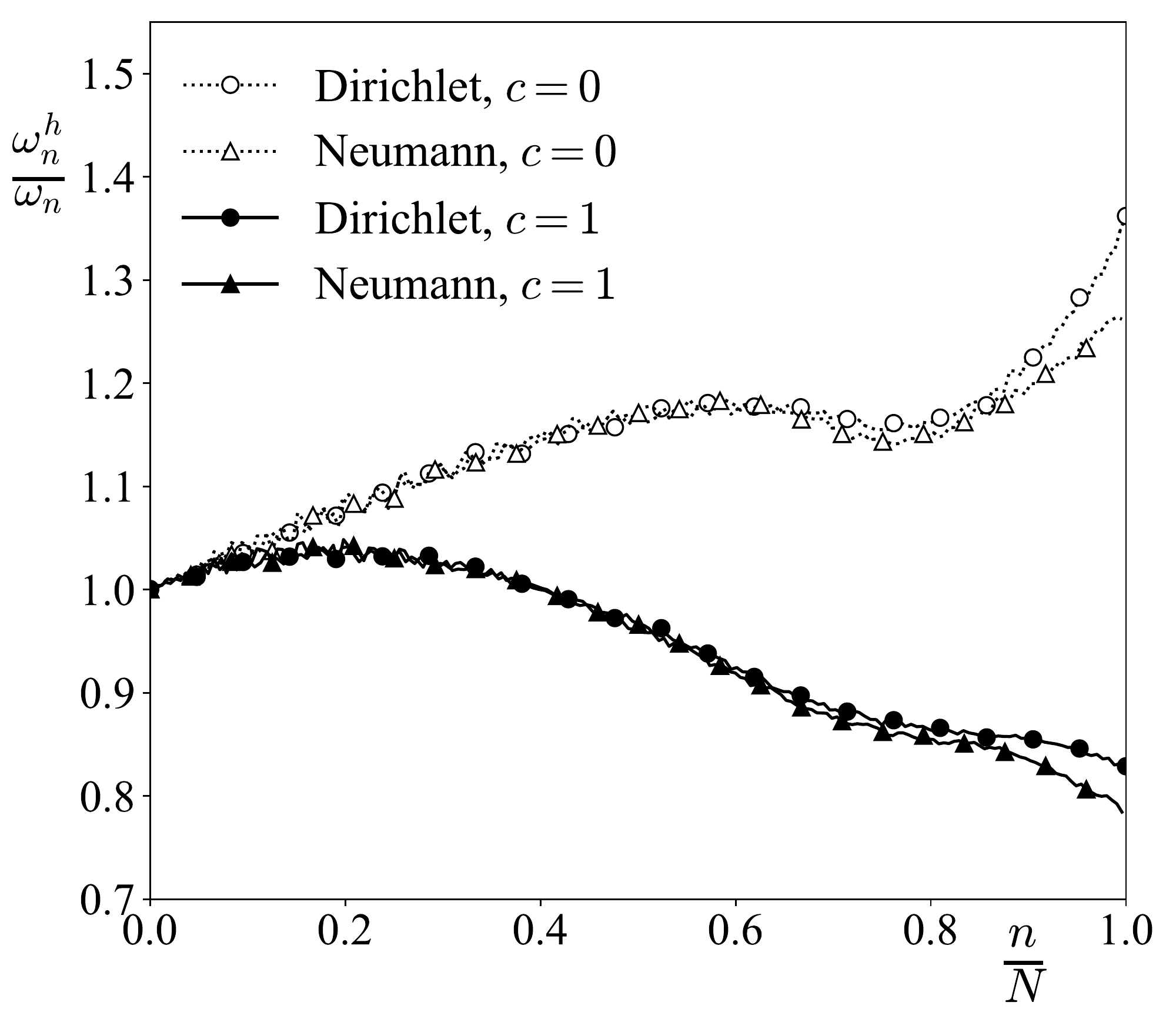}\label{2Dquads1}}\hfill
    \subfloat[$P=2$]{\includegraphics[width=0.5\linewidth,trim=0 13 0 10, clip]{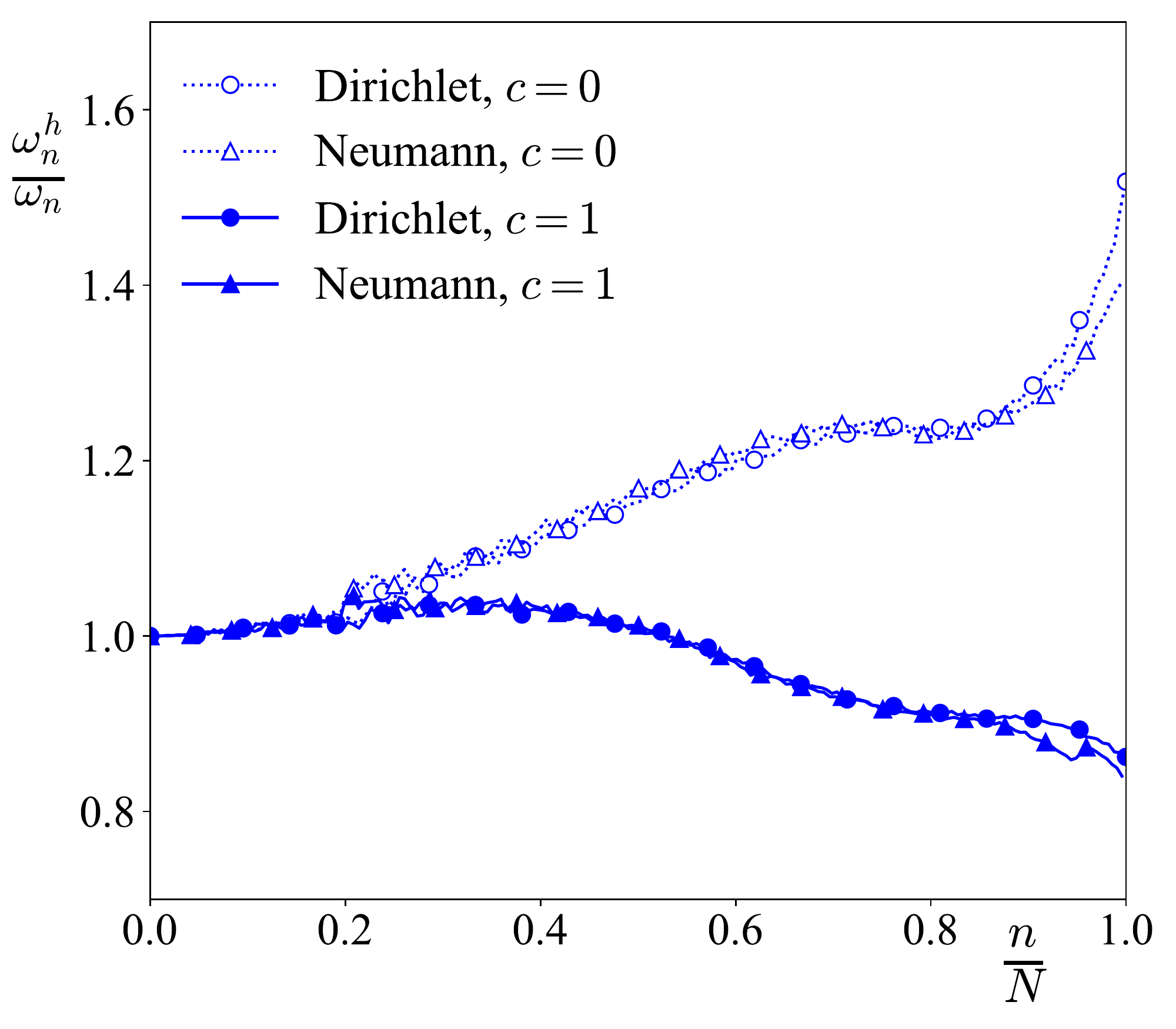}\label{2Dquads2}}\\
    \subfloat[$P=3$]{\includegraphics[width=0.5\linewidth,trim=0 13 0 10, clip]{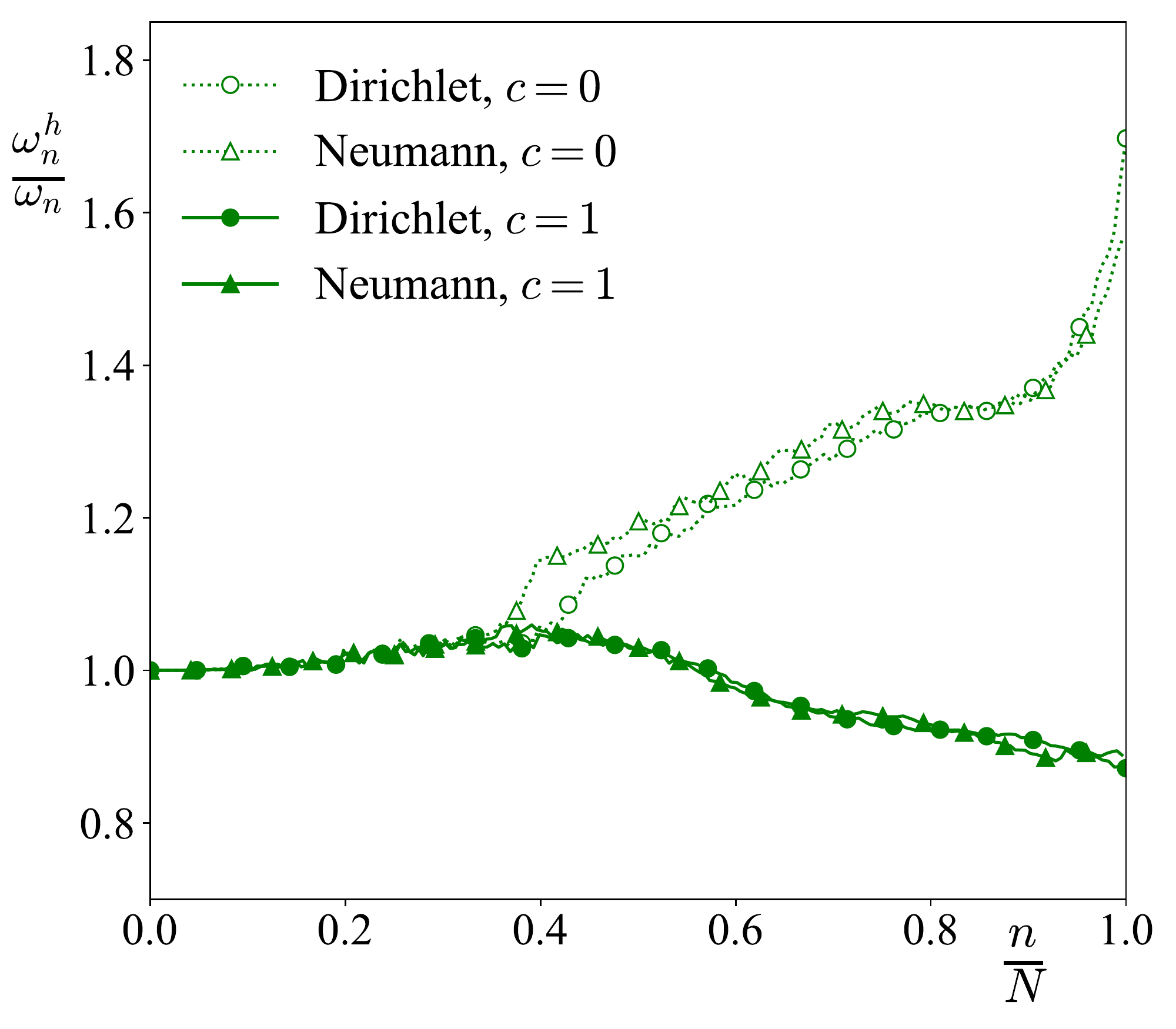}\label{2Dquads3}}\hfill
    \subfloat[$P=4$]{\includegraphics[width=0.5\linewidth,trim=0 13 0 10, clip]{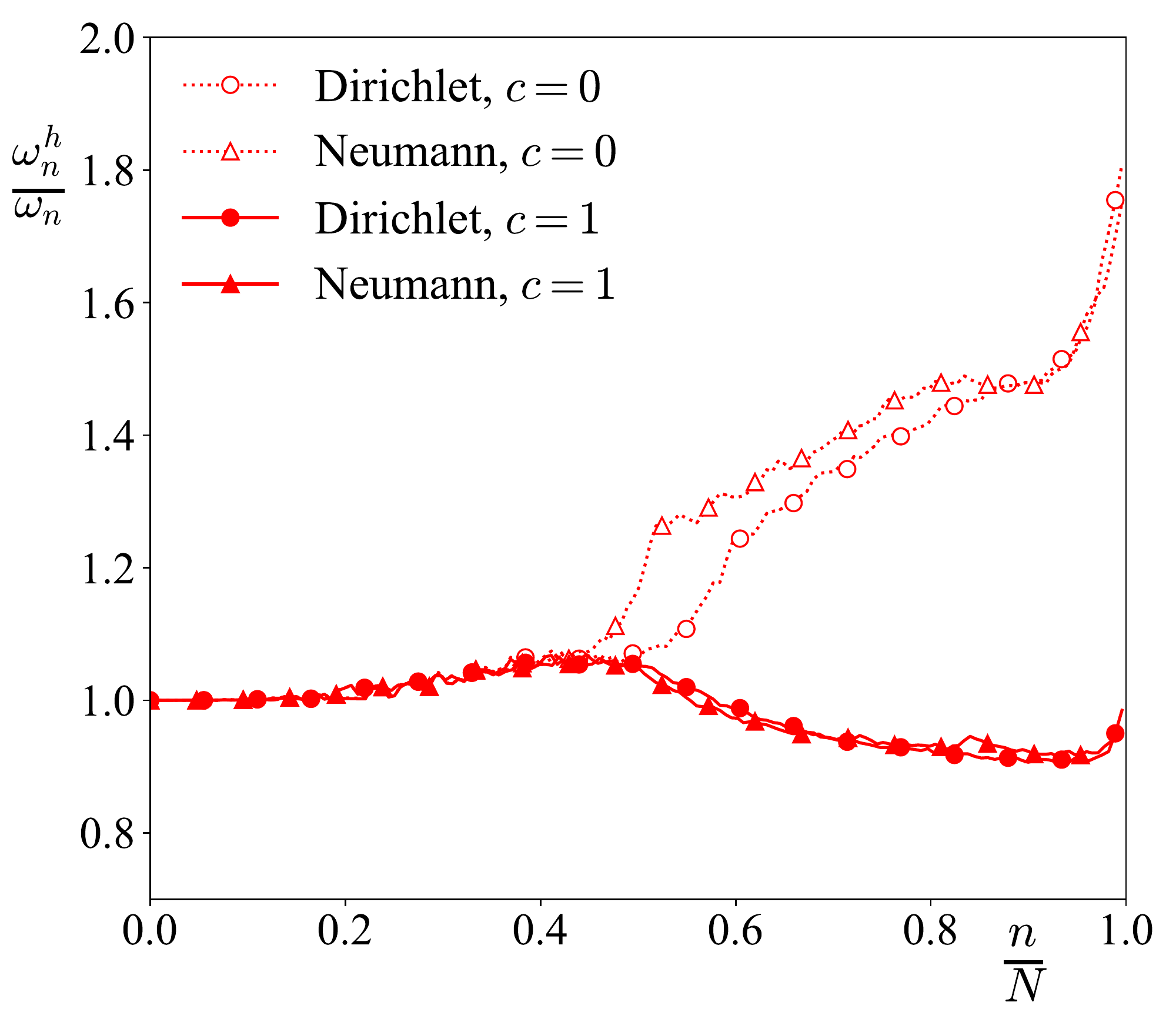}\label{2Dquads4}}
    \caption{Frequency spectra for a two-dimensional domain discretized with quadrilateral elements and various polynomial orders, normalized by the exact eigenfrequencies. All mass-scaled computations use~$c=1$.} \label{2Dquads}
\end{figure}

\begin{figure}[!b]\vspace{-0.25cm}
    \centering
    \subfloat[$P=1$]{\includegraphics[width=0.5\linewidth,trim=0 13 0 10, clip]{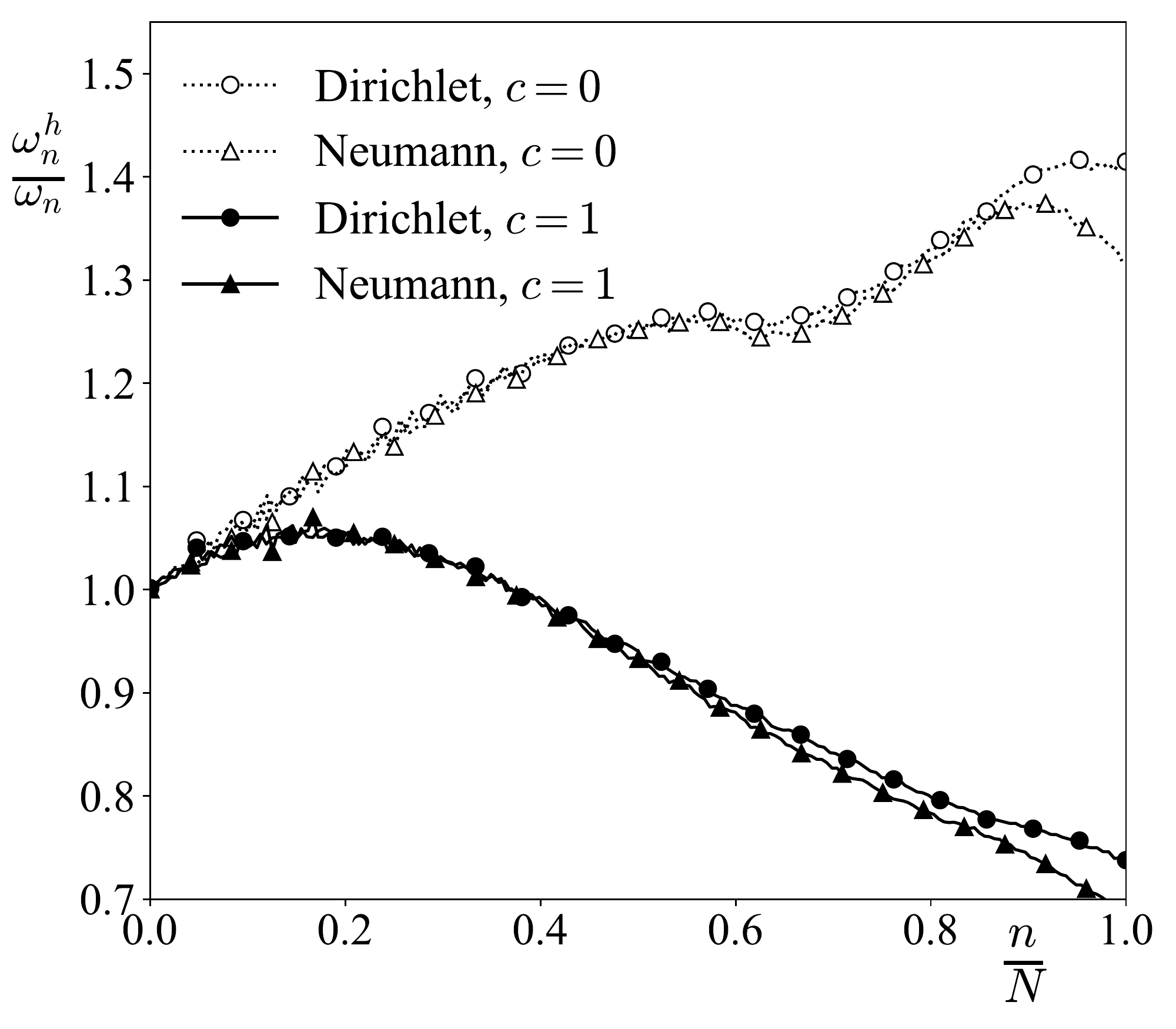}\label{2Dtriangles1}}\hfill
    \subfloat[$P=2$]{\includegraphics[width=0.5\linewidth,trim=0 13 0 10, clip]{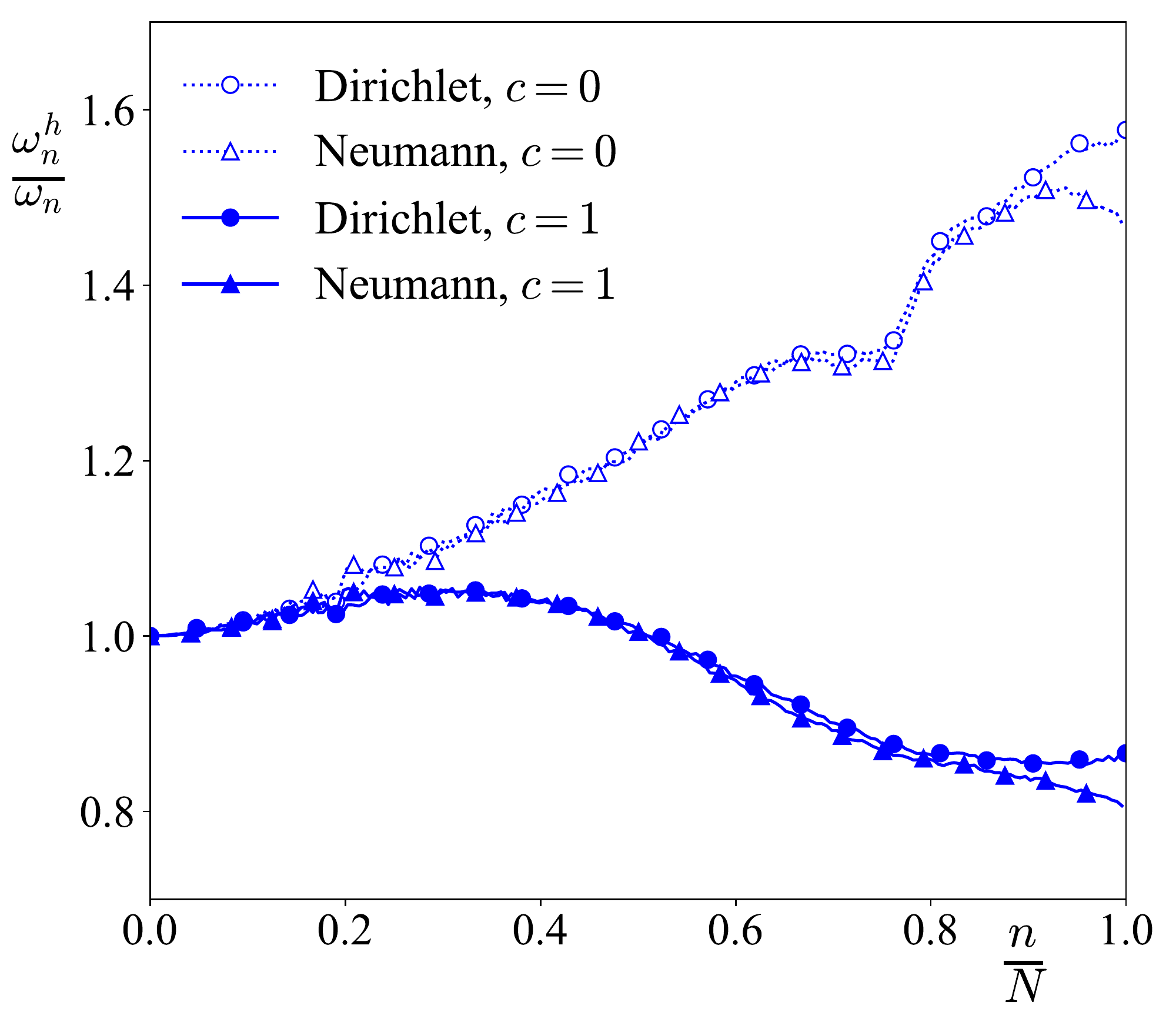}\label{2Dtriangles2}}\\
    \subfloat[$P=3$]{\includegraphics[width=0.5\linewidth,trim=0 13 0 10, clip]{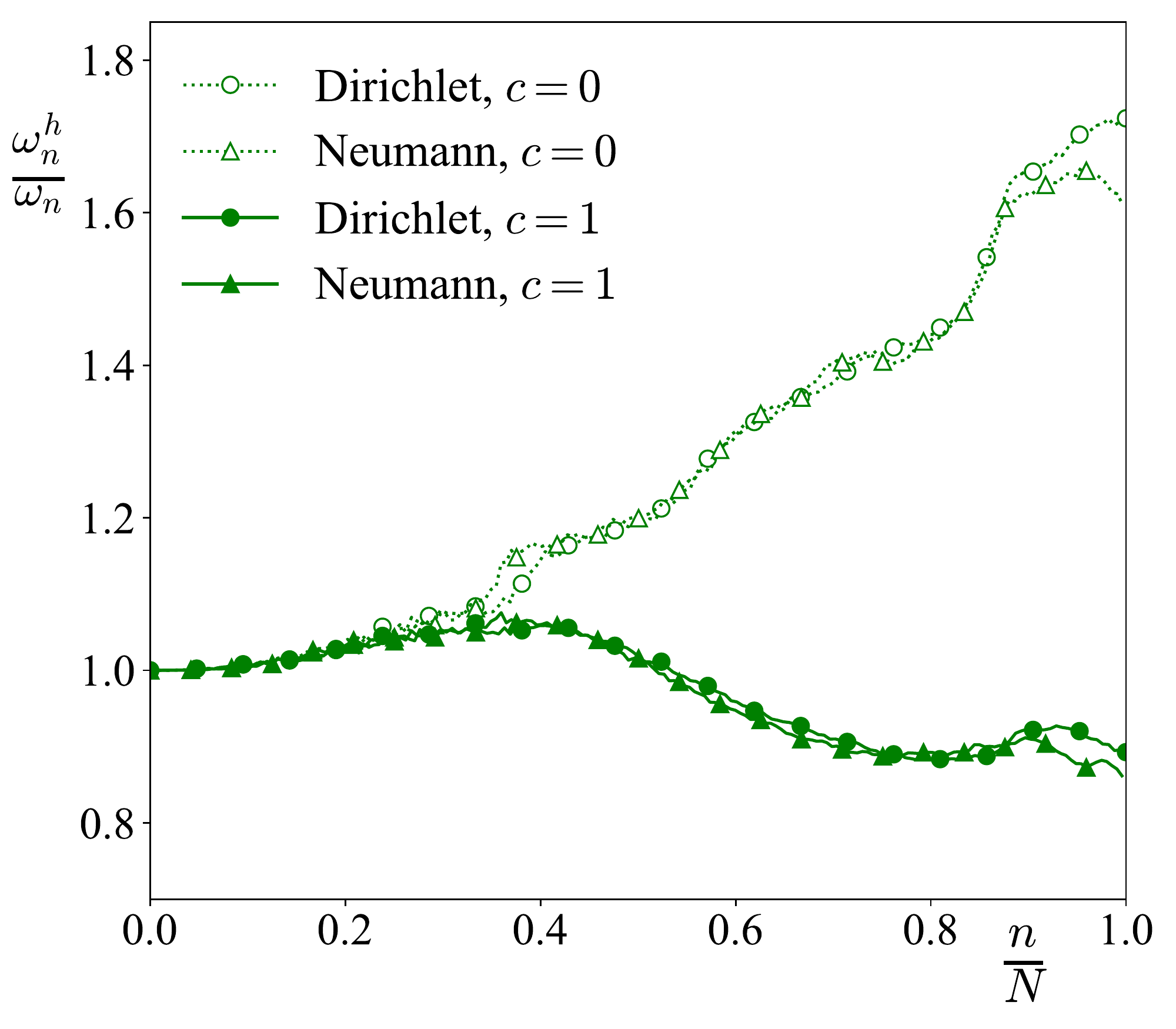}\label{2Dtriangles3}}\hfill
    \subfloat[$P=4$]{\includegraphics[width=0.5\linewidth,trim=0 13 0 10, clip]{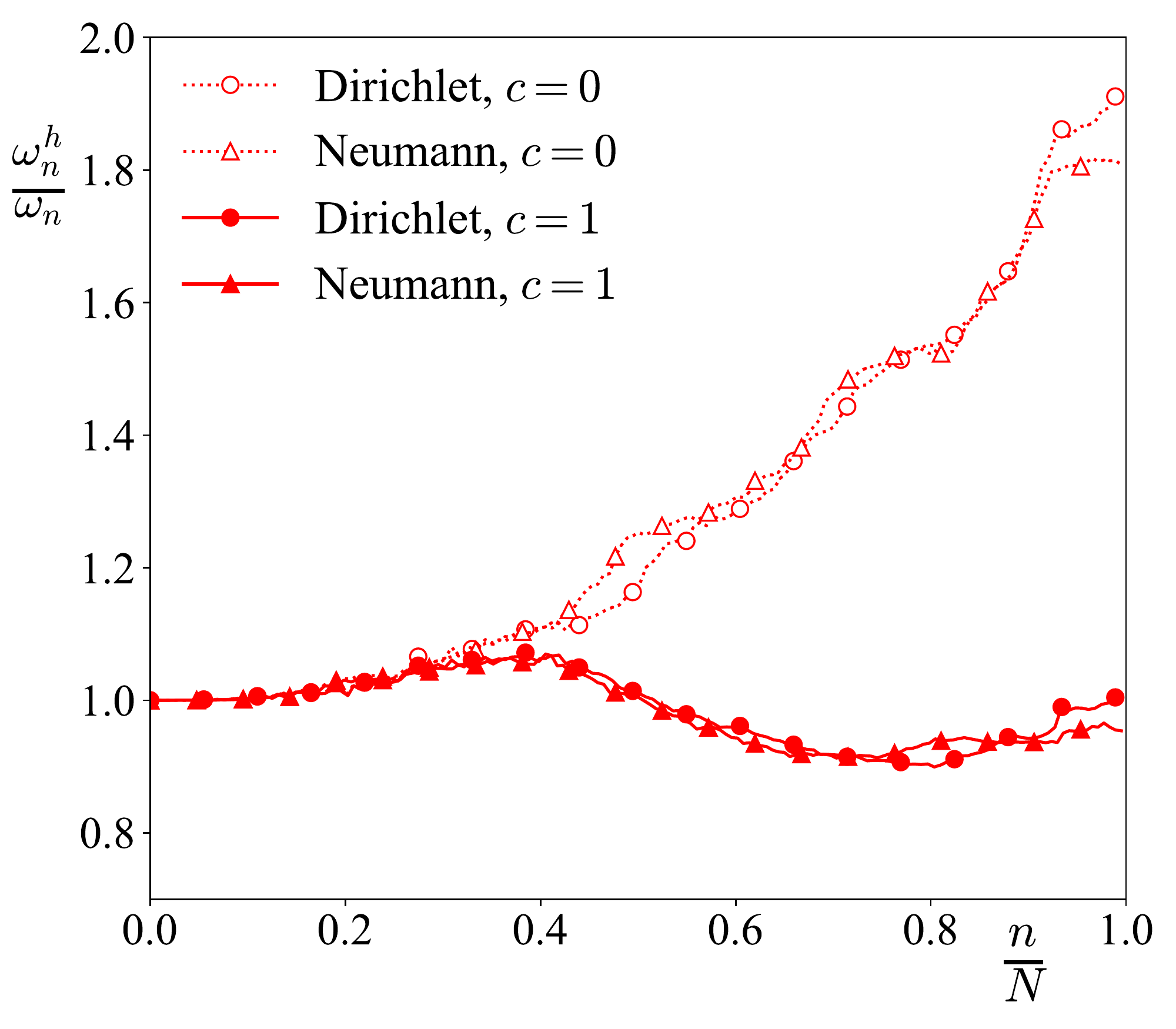}\label{2Dtriangles4}}
    \caption{Frequency spectra for a two-dimensional domain discretized with triangular elements and various polynomial orders, normalized by the exact eigenfrequencies.. All mass-scaled computations use $c=1$.}\label{2Dtriangles}
\end{figure}

Next, we consider a two-dimensional square domain with either Dirichlet or Neumann conditions around the entire circumference. The Neumann boundary case is particularly important, since the mass-scaling term of \cref{bilinpenalty} includes a component specifically on the Neumann boundary. Of course, for pure Neumann conditions there is a single zero eigenvalue corresponding to the constant solution eigenfunction. In the following, we remove this eigenvalue from the spectrum. 

\Cref{2Dquads1,2Dquads2,2Dquads3,2Dquads4} show the results for $P=1$ to $4$ when square elements are used and \cref{2Dtriangles1,2Dtriangles2,2Dtriangles3,2Dtriangles4} for triangular elements. In all cases, the numbers of elements are chosen such that the total number of degrees of freedom is kept at approximately 900. The overall behavior due to the mass scaling is very similar to that of the one-dimensional case, irrespective of the type of boundary condition and the type of elements used. We again observe that the lower end of the spectrum is unaffected, that accuracy of the frequencies in the medium range of the spectrum is significantly improved, and that the frequencies of the medium to higher end of the spectrum are decreased compared to the cases without mass scaling. The consistency between these results and the results for the one-dimensional case for various polynomial orders is remarkable, given the sensitivity to $c$ and the predicted non-trivial dependency on $P$, $d$ and $h$ in \cref{beta}.

%=====================================================
\section{Dynamic response: critical time step and convergence}
\label{sec:dynamics}
%=====================================================

\begin{figure}[!b]
    \centering
    \subfloat[Hexagonal cut-out.]{\includegraphics[width=0.46\linewidth]{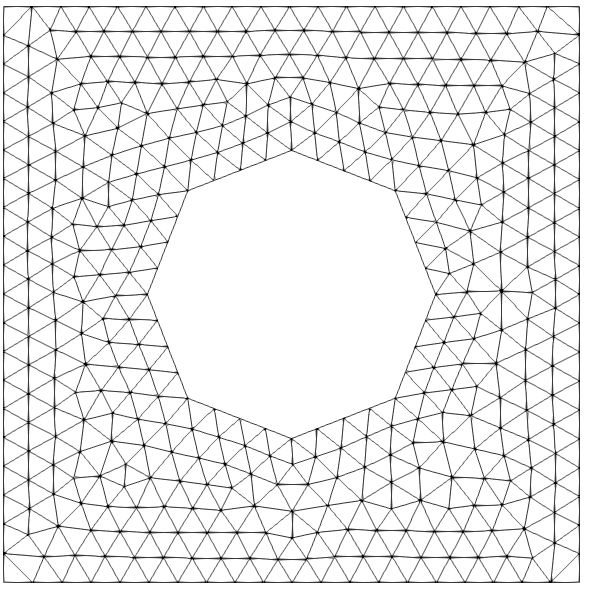} }\hspace{.25cm}
    \subfloat[Star-shaped cut-out with local refinement.]{\includegraphics[width=0.46\linewidth]{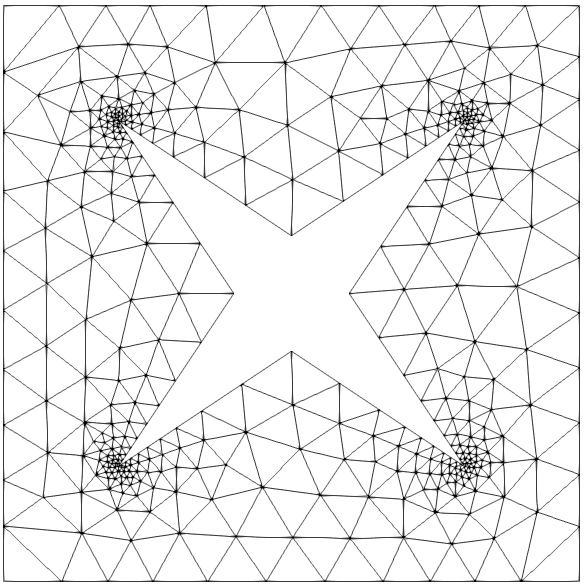} }
    \caption{Meshes with approximately $400$ vertices.}
     \label{fig:meshes}
\end{figure}

The examples of \cref{sec:spectra} concern tensor-product domains and uniform meshes with regularly shaped elements. To understand if the improvement of the spectrum holds practical significance for more complex cases, we consider a square domain with a polygonal cut-out, and compute the dynamic response of a forced vibration. The improvement of the spectrum should permit a larger time step when explicit time integration is employed.
To understand the generality of the impact on the critical time step, we consider two types of meshes with different cut-outs: a mesh with an octagonal cut-out and uniformly distributed elements and a mesh with a star-shaped cut-out with local refinement at the star points. Both meshes are illustrated in \cref{fig:meshes} for the refinement case with approximately 400 vertices. 

\begin{figure}[!b]
    \centering
    \subfloat[$t_{\text{crit}}$.]{\includegraphics[width=0.49\linewidth]{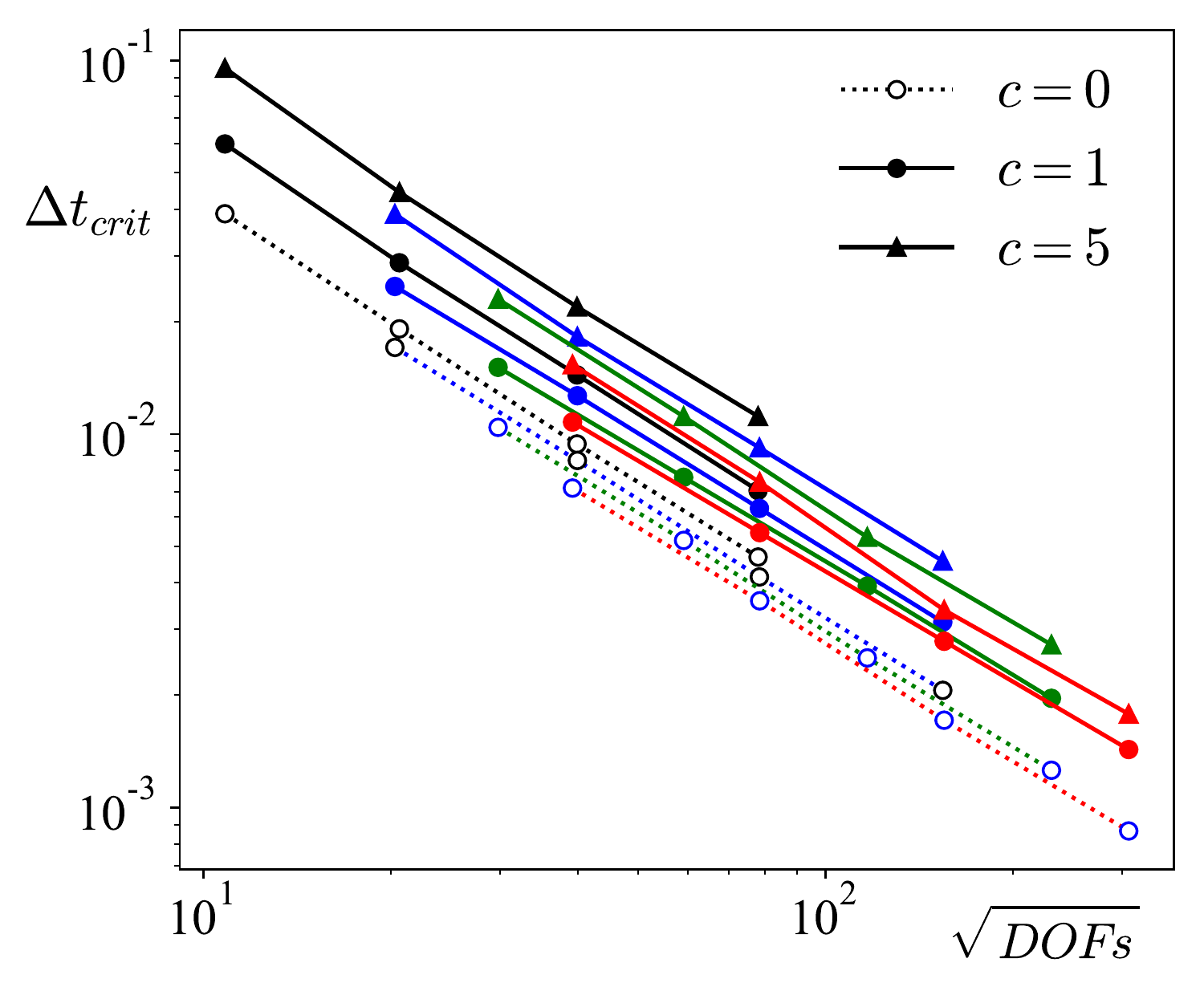} }\hfill
    \subfloat[Increase in $t_{\text{crit}}$.]{\includegraphics[width=0.49\linewidth]{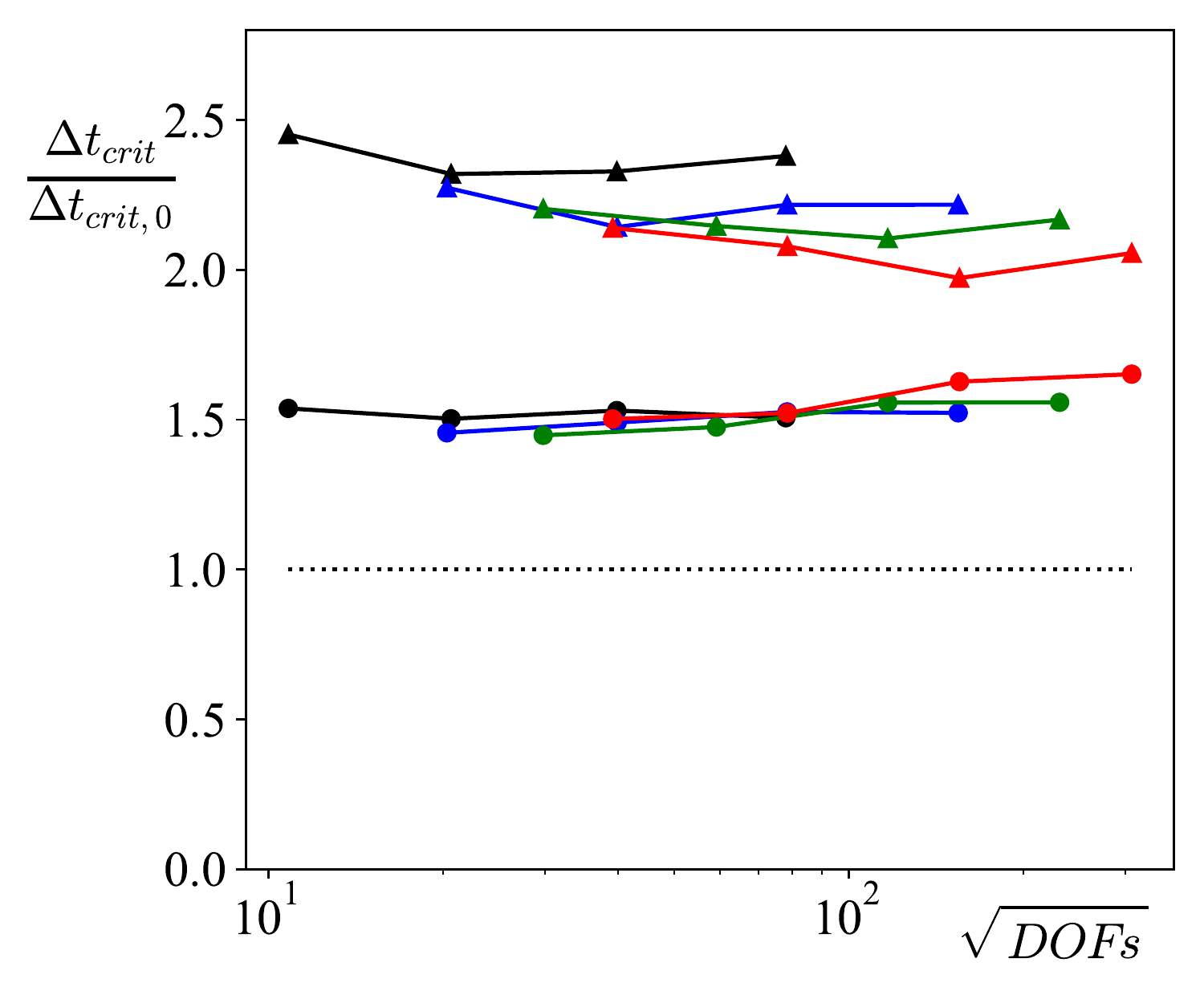} }
    \caption{Critical time step for a square domain with an octagonal cut-out.}\label{fig:ct_oct}
\end{figure}

\begin{figure}[!b]
    \centering
    \subfloat[$t_{\text{crit}}$.]{\includegraphics[width=0.49\linewidth]{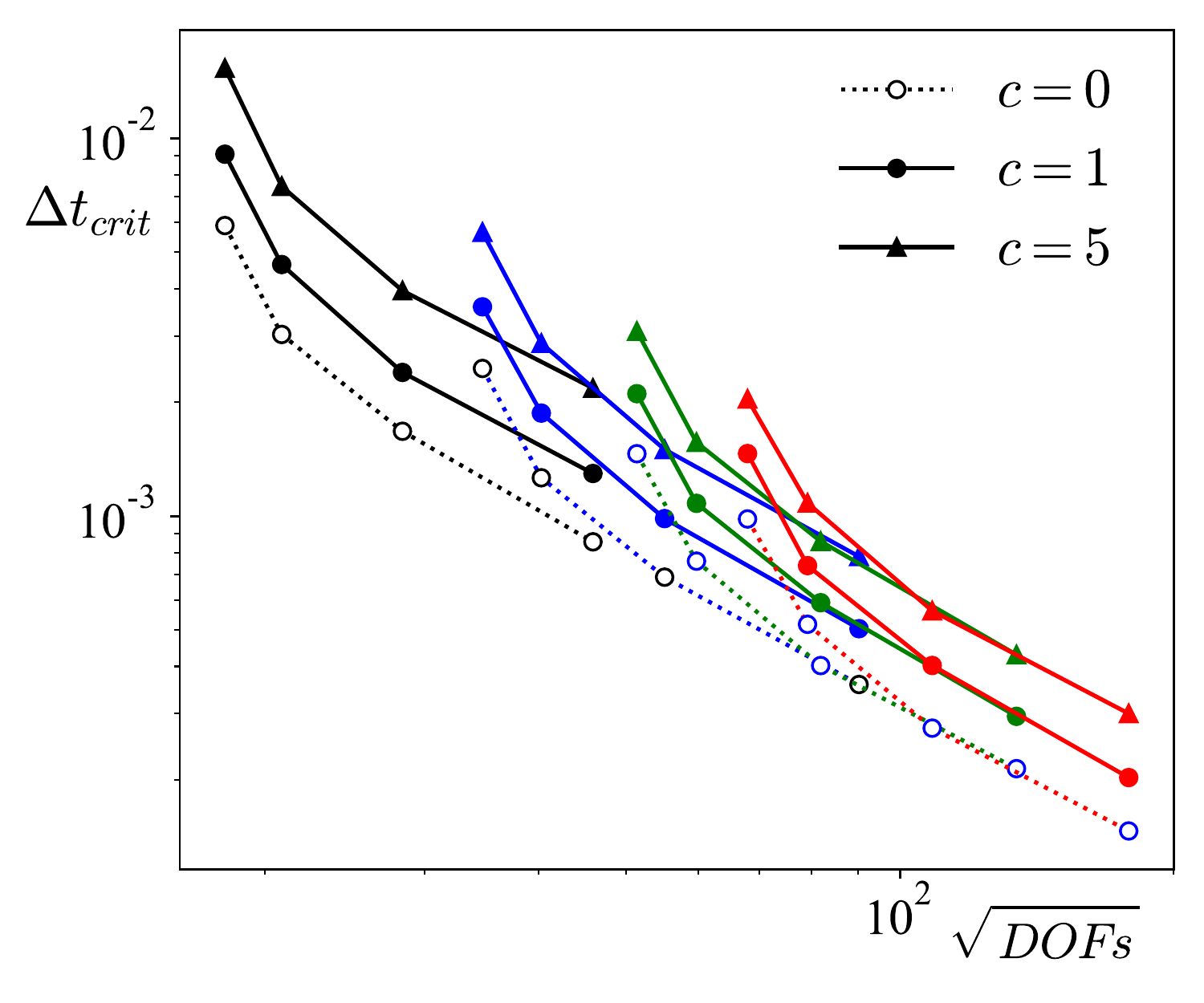} }\hfill
    \subfloat[Increase in $t_{\text{crit}}$.]{\includegraphics[width=0.49\linewidth]{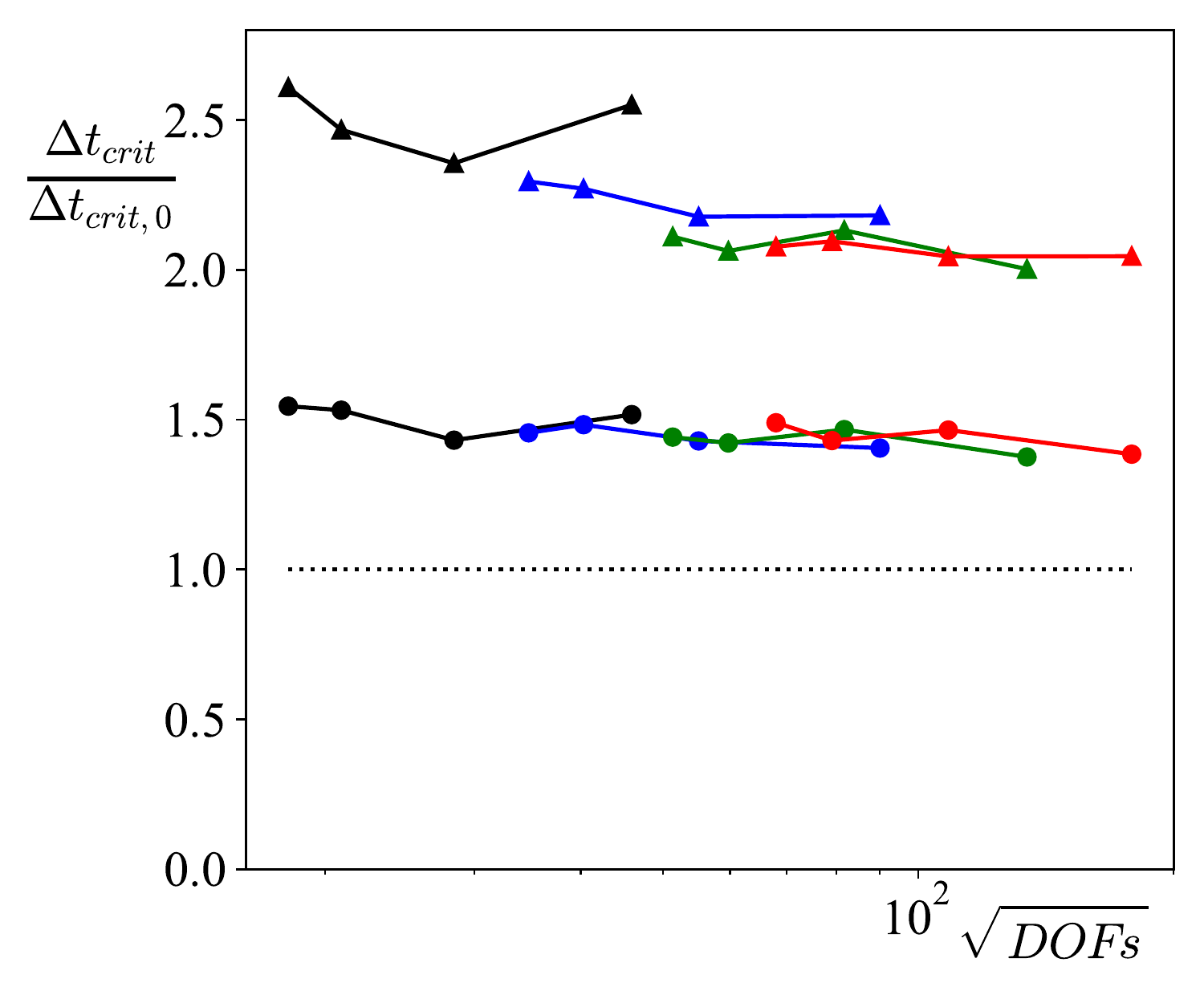} }
    \caption{Critical time step for a square domain with a star-shaped cut-out and local mesh refinement.}\label{fig:ct_star}
\end{figure}

We compute the critical time step with \cref{Dtcrit} for multiple refinement-levels and polynomial orders. The results are shown in \cref{fig:ct_oct,fig:ct_star} for the octagonal cut-out and the star-shaped cut-out, respectively. The left figures show the values of the critical time steps and the right figures show the factor by which they increase due to mass scaling. The black, blue, green and red lines indicate polynomial orders of 1, 2, 3 and 4, respectively, and the dotted line with the open markers concerns the case without mass scaling (i.e., $c=0$). When mass scaling is added, we observe a uniform increase in critical time-step size for both domains, all refinements and all polynomial orders. With a larger value of $c$, the increase in time-step size is larger. The increase by $50\%$ for $c=1$ corresponds well to the prediction from \cref{ssec:abc}. The obtained increase by $100\%-150\%$ for $c=5$ is lower than the predicted value and the prediction becomes less accurate for the higher polynomial orders. We anticipate that this overprediction is related to the reordering of the maximum eigenvalues. This was discussed in \cref{ssec:1D} and observed in \cref{fig:1D_abs}. An increase of $c$ suppresses the eigenvalues of the highest frequency eigenfunction most, which, at a certain $c$, causes the eigenvalue corresponding to a lower frequency eigenfunction to dominate. This reduces the effectiveness of the mass-scaling term. As a result, the approximation strategy of \cref{ssec:abc}, which may be considered a first-order approximation around $c=0$, becomes an overprediction.

Next, we compute the dynamic response for a system with the following true solution:
\begin{align}
    \begin{cases}
    & u(t,x,y) = \sum\limits_{i=1}^{9} w_i \cos\left( \sqrt{n_i^2+m_i^2} \,\pi\, t \right) \sin\left(n_i \pi x \right) \sin\left( m_i \pi y \right)  \,, \\
    & w = \{ 1,0.8,0.8,0.6,0.5,0.2,0.1,0.05,0.03 \}  \,, \\
    & n = \{ 3,4,6,8,5,12,9,12,8 \}  \,, \\
    & m = \{ 4,3,8,6,12,5,12,9,15 \} \,, 
    \end{cases}
\end{align}
from which we determine the initial condition, as well as the Neumann data $g$ and $\ddot{g}$ on the circumference of the cut-out. This true solution has a period of 2 seconds.

\begin{figure}[!b]
    \centering
    \subfloat[$T_{\text{max}} = 0.1$ and $\Delta t = 0.9\,t_{\text{crit}}$.] {\includegraphics[width=0.49\linewidth]{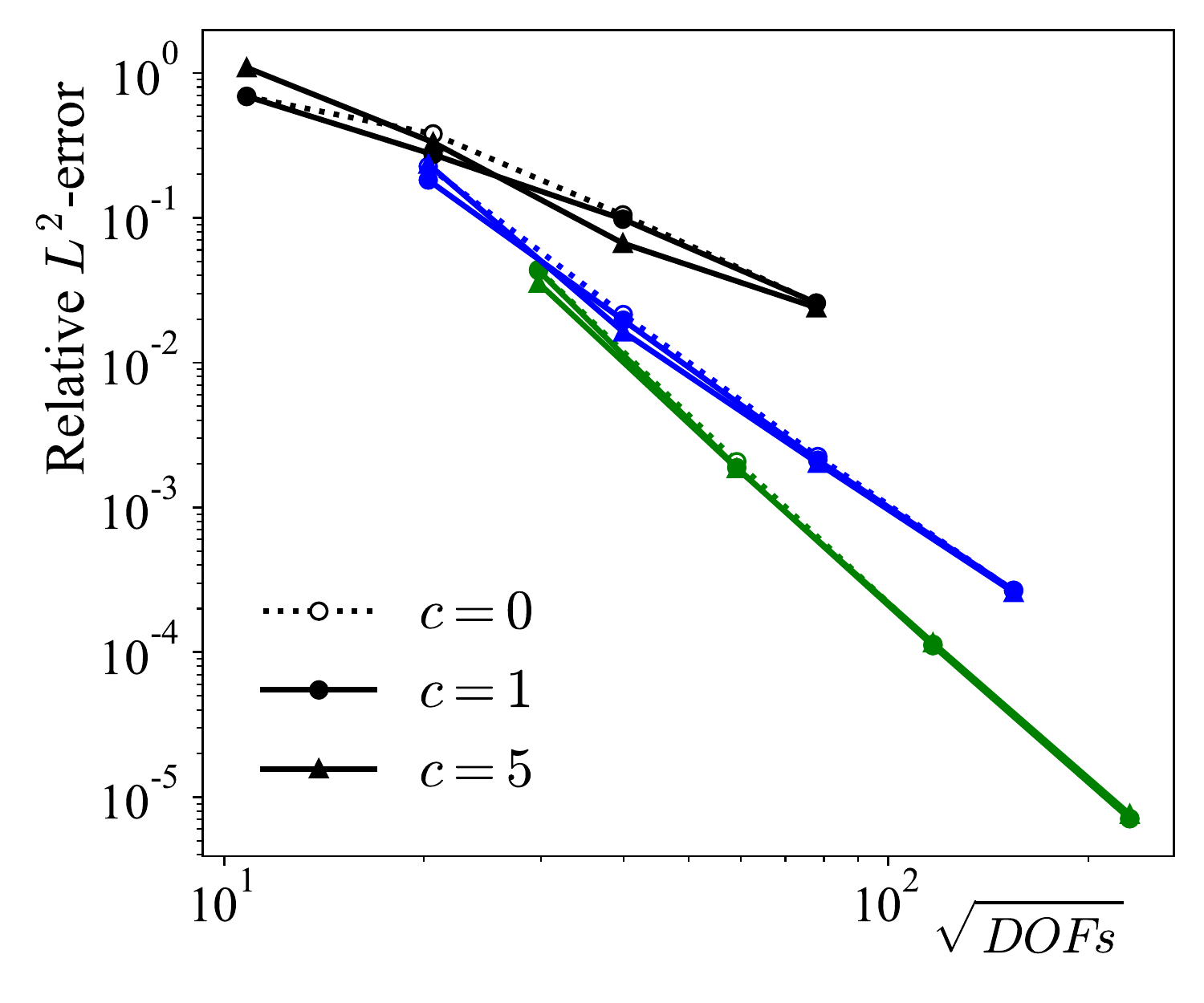} \label{fig:conva}}\hfill
    \subfloat[$T_{\text{max}} = 10.1$ and $\Delta t = 0.9\,t_{\text{crit}}$.] {\includegraphics[width=0.49\linewidth]{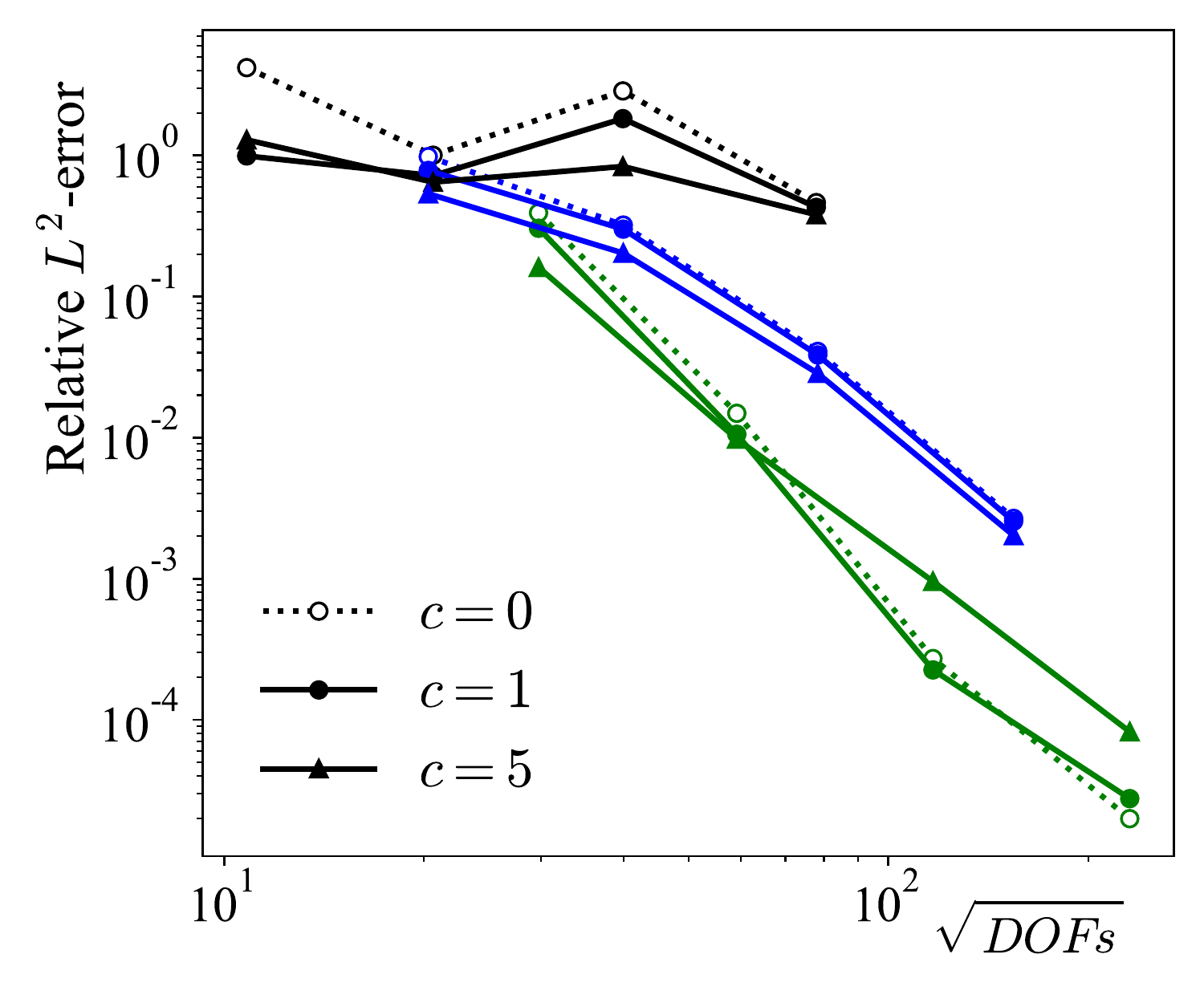} \label{fig:convb}}
    \caption{$L^2$-convergence for the square domain with an octagonal cut-out.}
\end{figure}

For the convergence study we focus on the case of the octagonal cut-out, as the refinement levels involve uniform refinement. We make use of the standard fourth-order explicit Runge-Kutta algorithm for the time integration. \Cref{fig:conva} shows the convergence of the $L^2$-error after $T=0.1$ seconds and \cref{fig:convb} shows the results precisely five periods later, at $T=10.1$ seconds. The previously plotted critical time steps are used as time steps in the Runge-Kutta algorithm. This means that the data points for $c=1$ and $c=5$ required approximately 33\% and 60\% fewer time step computations than those without mass scaling. In \cref{fig:conva} we observe almost identical errors and optimal rates of convergence for all cases, indicating that the proposed mass scaling does not negatively affect the solution quality. The convergence curves plotted in \cref{fig:convb} are less clean. Nevertheless, the mass scaling never negatively affects the errors, despite the increase in time-step size. The exceptions are the last two datapoints for the case of $P=3$ and $c=5$, which could be due to the error in the time discretization dominating the complete solution error at this level of spatial resolution.

%=====================================================
\section{Extension towards non-linear solid mechanics}
\label{sec:3Dexample}
%=====================================================

To illustrate the use of our method for a more complex scenario, we consider a beam oscillating at high amplitudes compared to its length. The non-negligible rotational deformation invalidates the linear strain assumption. As a consequence, a geometrically non-linear model needs to be considered. In the domain $\Omega$ of the reference configuration, the variational statement governing free vibration reads \cite{belytschko2014,Zienkiewicz:00.2}:
\begin{align}
&\text{Find }\B{u}\in \B{H}^1\big(\mathcal{T},\B{H}^1_{\boldsymbol{0}}(\Omega)\big) \text{ s.t., for a.e. }t,\,\forall\,\B{v} \in \B{H}^1_{\boldsymbol{0}}(\Omega):\nonumber\\
&    \quad \big\langle \rho_0\frac{\partial^2 }{\partial t^2}\B{u} ,  \B{v} \big\rangle_{\boldsymbol{H}^{-1},\boldsymbol{H}^{1}_0}  + \int\limits_\Omega \nabla \B{v} : \B{P}(\B{u}) \dO = 0  \,, \label{NLEl}
\end{align}
where $\B{u}$ is the displacement vector, $\rho_0$ the density in the reference configuration and $\B{P}(\B{u})$ the first Piola-Kirchoff stress tensor. The first and second Piola-Kirchoff stress tensors relate to one-another via the deformation gradient:
\begin{align}
    \B{P}(\B{u}) = \B{F}(\B{u})\cdot\B{S}(\B{u}) = ( \B{I} + \nabla\B{u}) \cdot\B{S}(\B{u})  \,.  \label{PK1PK2}
\end{align}
As a constitutive law we make use of the Saint Venant–Kirchhoff hyperelastic relation:
\begin{align}
   \B{S}(\B{u})  = {\mathbb{C}}:\B{E}(\B{u}) =   \lambda \,\text{tr}(\B{E}(\B{u}))\B{I} + 2\mu\B{E}(\B{u}) \,, 
\end{align}
with $\mathbb{C}$ the fourth order material response tensor, $\lambda$ and $\mu$ the Lam\'e parameters and $\B{E}(\B{u})$ the Green-Lagrange strain tensor:
\begin{align}
\B{E}(\B{u}) = \frac{1}{2} \big( \B{F}(\B{u})^\text{T}\cdot\B{F}(\B{u})  - \B{I} \big) =  \frac{1}{2}\big( \nabla \B{u} + \nabla \B{u}^\text{T} + \nabla \B{u}^\text{T}\cdot \nabla \B{u} )  \,. \label{EGL}
\end{align}
All gradients in \cref{NLEl,PK1PK2,EGL} operate in the reference configuration.

Under the assumption that the solution is smooth enough we may perform integration by parts on the second term to reveal the implicitly assumed transmission conditions between two neighboring domains:
\begin{align}
    \jump{ \B{P}(\B{u}) \B{n} } = \B{0} \qquad \text{on }\Gamma  \,, \label{PKjump}
\end{align}
so that it is the continuity of the first Piola-Kirchoff stress that represents the interfacial equilibrium condition when observed in the reference configuration \cite[p.~204]{belytschko2014}. 

The theory outlined in \cref{sec:formulation} thus suggests to use a discrete formulation of \cref{NLEl} that also includes a mass-scaling term based on the second time derivatives of \cref{PKjump}. However, due to its non-linearity, addition of a penalty of this traction jump to the mass matrix would require re-assembly and re-computation of the LU-decomposition of the stabilized mass matrix in each time step. The required additional computational expense to perform these operations is of same order of magnitude as the expense of a single implicit time-integration step, and this thus defeats the point of using an explicit time-integration method. Instead, we propose to add a linearized version of this penalty. Then, the extension of \cref{FEM} to the current problem would produce:
\begin{align}
&\text{Find }\B{u}^h\in \B{U}^h \text{ s.t. }\forall\,\B{v}^h \in \B{U}^h: \nonumber\\
\begin{split}
&    \quad \int\limits_\Omega \rho_0 \frac{\partial^2 }{\partial t^2} \B{u}^h \cdot \B{v}^h \dO + \int\limits_\Gamma \beta \frac{ \rho }{ |\mathbb{C}|^2 } \jump{\big(\mathbb{C}: \boldsymbol{\epsilon}(\tfrac{\partial^2 }{\partial t^2}\B{u}^h)\big) \B{n} }\cdot \jump{ \big(\mathbb{C}: \boldsymbol{\epsilon}(\B{v}^h)\big) \B{n} }\dG \\
& \hspace{1cm}+\int\limits_{\partial\Omega_N} \beta \frac{ \rho }{ |\mathbb{C}|^2 } \Big[\big(\mathbb{C}: \boldsymbol{\epsilon}(\tfrac{\partial^2 }{\partial t^2}\B{u}^h)\big) \B{n} \Big]\cdot \Big[ \big(\mathbb{C}: \boldsymbol{\epsilon}(\B{v}^h)\big) \B{n} \Big]\dG  + \int\limits_\Omega \nabla \B{v}^h : \B{P}(\B{u}^h) \dO = 0  \,, 
\end{split} \label{formulationNLEL}
\end{align}
where we have assumed homogeneous Dirichlet and Neumann conditions. For $|\mathbb{C}|$ we take its maximum eigenvalue, that is $\max(\mu,K)$ with $K$ the bulk modulus as $K=\lambda+\frac{2}{3}\mu$. For $\beta$, we make use of a straightforward extension of the expression from \cref{beta} from the scalar wave equation to a vector wave equation:
\begin{align}
\beta = c \frac{1}{4}\frac{1}{d^4 \pi^2 }\frac{1}{2P^3-P^2} h^3 \,. \label{beta_vec}
\end{align}
The additional division by $d^2$ is due to each node carrying $d$-degrees of freedom, such that the index of the maximum eigenvalue as used in \cref{ssec:trueeig} now becomes $(dPn)^d$. 
%which means that the maximum eigenvalues in \cref{eigenvalues_num} corresponds to the $(dPn)^d$-th eigenvalues.

The new terms that are added to the mass matrix involve the linear strain:
\begin{align}
    \boldsymbol{\epsilon}(\B{u}) = \frac{1}{2}\big( \nabla \B{u} + \nabla \B{u}^\text{T} \big) \,. 
\end{align}
This raises questions regarding the variational consistency of the formulation. The linearized strain measure yields erroneous (i.e., nonzero) strain values for rigid body rotation. However, this is not necessarily problematic: we add a penalty to the \textit{jump} of this linearized stress, rather than the linearized stress itself. When the material parameters are continuous, and the true solution is smooth enough, \cref{formulationNLEL} is still satisfied by the analytical solution. We investigate the impact of this model simplification by focusing on an example that involves non-negligible rotations.

\begin{figure}[!b]
\vspace{-0.5cm}
    \centering
    \subfloat[Beam geometry.]{\includegraphics[width=0.49\linewidth]{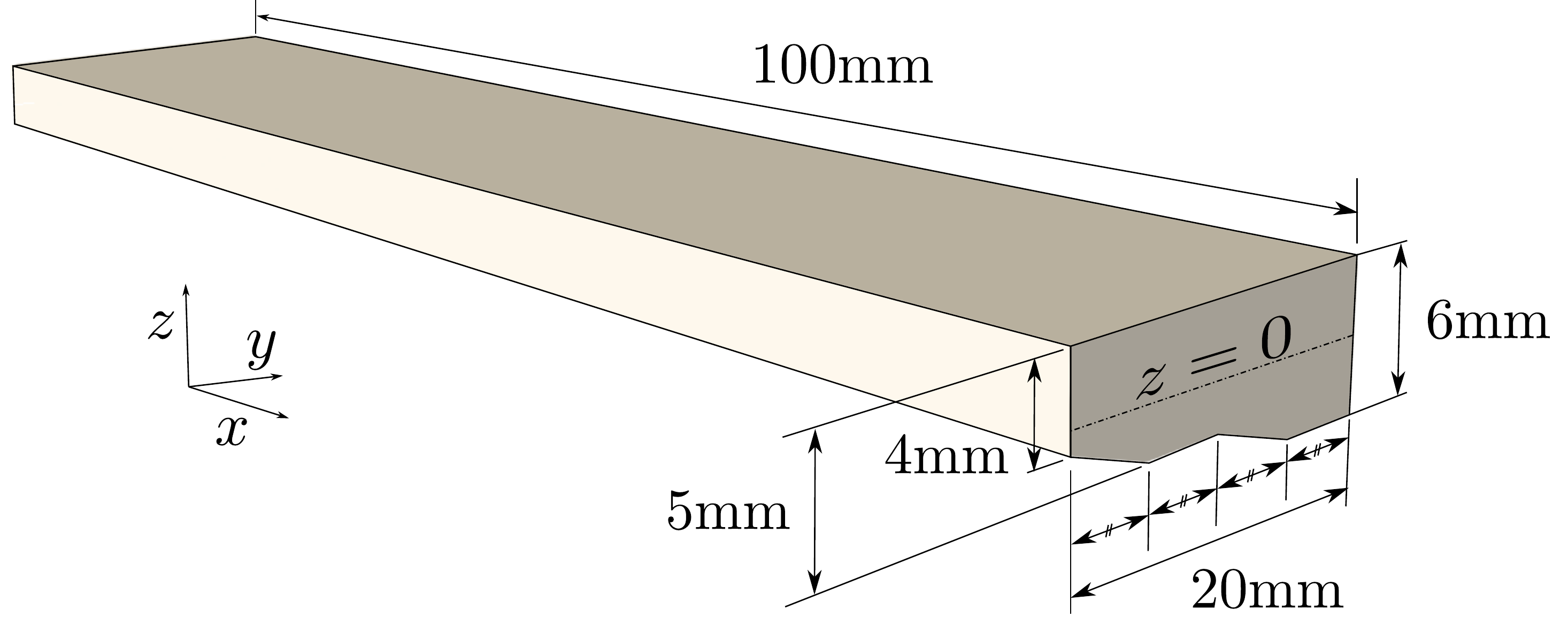} \label{fig:beam} }
    \subfloat[Initial displacement.]{\includegraphics[width=0.45\linewidth]{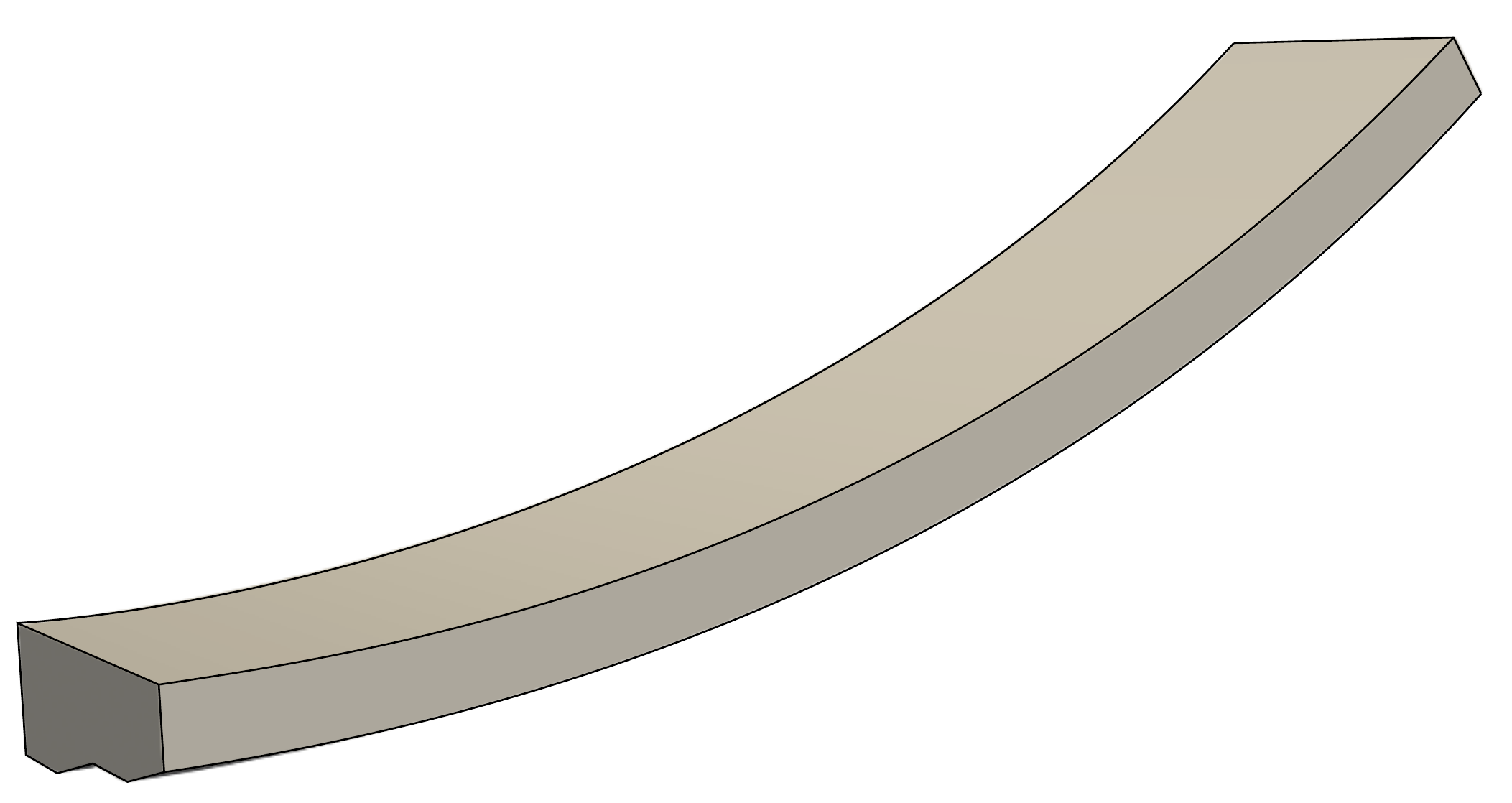} \label{fig:beamIC} }\\
    \subfloat[Coarse mesh.]{\includegraphics[width=0.49\linewidth]{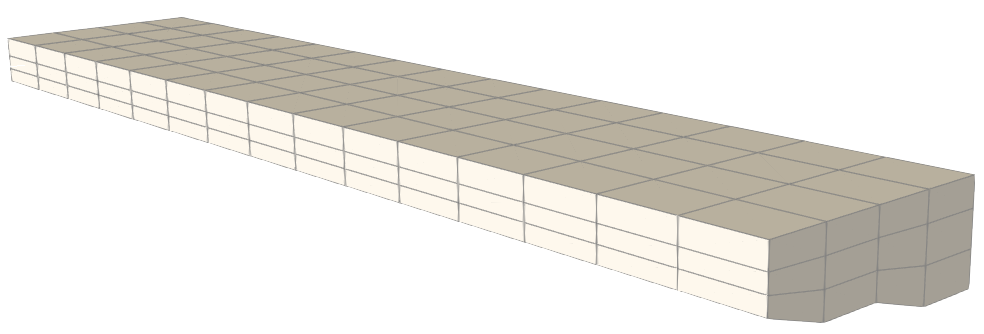} \label{fig:3Dmesh1} }
    \subfloat[Fine mesh.]{\includegraphics[width=0.49\linewidth]{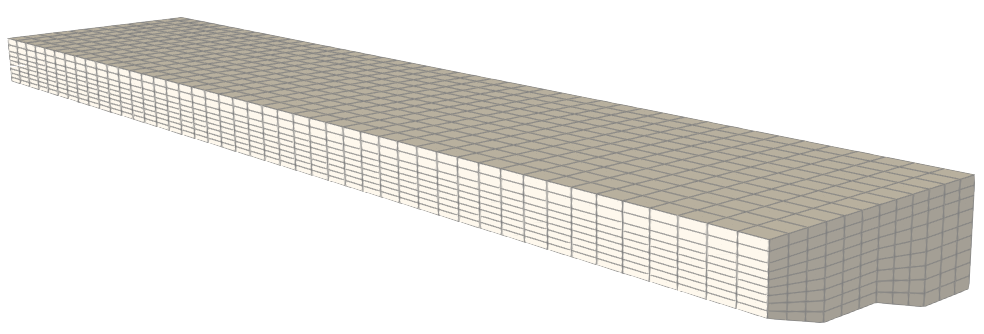} \label{fig:3Dmesh2} }
    \caption{Three-dimensional model problem.}\label{fig:3Dgeom}
\end{figure}

Consider the asymmetric beam illustrated in \cref{fig:beam}. The material parameters are chosen as $\rho=2700$kg/m$^3$, $E=73$GPa/m$^2$ and $\nu=0.33$, and correspond to an aluminum alloy. The beam is clamped at $x=0$, which is the left-most face in the illustrations of \cref{fig:3Dgeom}. At $t=0$, the beam is released from the following stationary initial displacement:
\begin{align}
    \B{u}_{\text{IC}} = \begin{bmatrix} (R-z)\sin(x/R)-x \\ 0 \\R-(R-z)\cos(x/R)-z \end{bmatrix} \,, 
\end{align}
for $R=0.15$. This describes a circular shape in the $y=0$ plane with a radius $R$ with no stretching in the initial $z=0$ plane. This initial displacement is shown in \cref{fig:beamIC}.
\Cref{fig:3Dmesh1,fig:3Dmesh2} illustrate the coarse and fine meshes used in the following computations. 

\begin{figure}[!t]
\vspace{-0.5cm}
    \centering
    \includegraphics[trim=10 15 10 10, clip,width=0.6\textwidth]{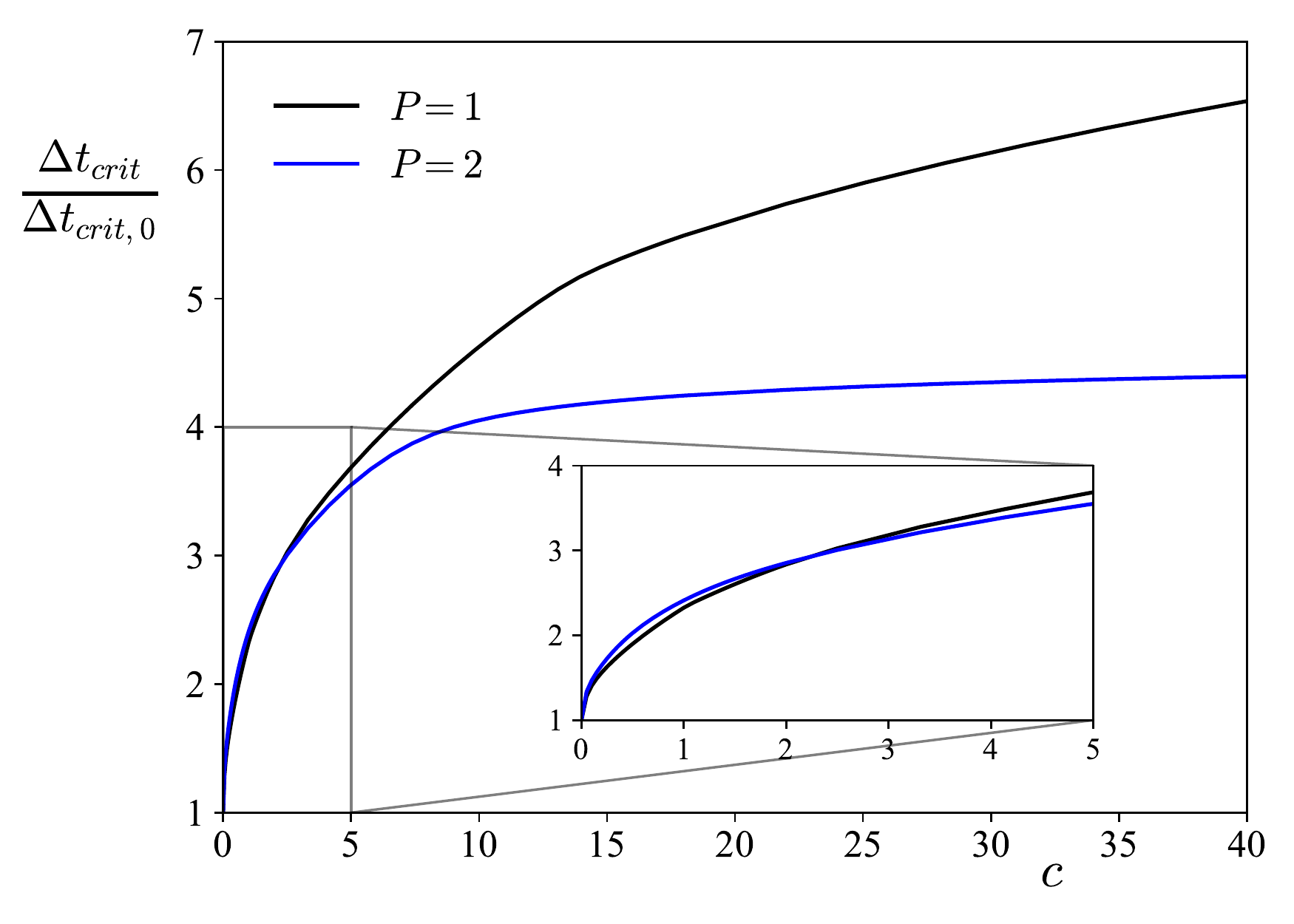}
    \caption{Factor of increase of the critical time step.\vspace{-0.5cm}}
    \label{fig:3Dtcrit}
\end{figure}

First, we study the dependency of the critical time step on $c$, the remaining free parameter in the estimate of $\beta$ in \cref{beta_vec}. To determine the eigenvalues for this nonlinear problem, we consider the stiffness matrix that arises from the linearization by a small perturbation around the initial displacement:
\begin{subequations}
\begin{alignat}{2}
    &K_{ij} = B(\B{N}_i,\B{N}_j) \,,\\
    \begin{split}
        &B(\B{v},\B{u}) = \int\limits_\Omega \nabla \B{v} : \big[ \nabla \B{u} \cdot \mathbb{C} : \B{E}(\B{u}_\text{IC})\big] \dO + \int\limits_\Omega \nabla \B{v} :  \big[ \B{F}( \B{u}_\text{IC} ) \cdot \mathbb{C} : \boldsymbol{\epsilon}(\B{u})\big] \dO\\
        &\hspace{2.5cm}+\frac{1}{2}\int\limits_\Omega  \nabla \B{v} :  \big[ \B{F}( \B{u}_\text{IC} ) \cdot \mathbb{C} : \big(\nabla \B{u}_\text{IC}^\text{T}\cdot \nabla \B{u} + \nabla \B{u}^\text{T} \cdot \nabla \B{u}_\text{IC}\big)\big] \dO \,.
    \end{split}
\end{alignat}
\end{subequations}
\Cref{fig:3Dtcrit} then shows the ratio between the critical time step and the critical time step for the case without mass scaling (i.e., $c=0$) as a function of $c$. These are computed for the coarse mesh. Overall, we observe significantly higher values for the increase in critical time step for this benchmark problem than for the benchmark problems of \cref{sec:spectra,sec:dynamics}. %This may be due to an inappropriate scaling of $\beta$. Recall that the expression of \cref{beta} is based on eigenvalue estimation for the scalar wave equation. 
The ratio keeps increasing when linear basis functions are used, whereas it levels off for the quadratic basis functions. This corresponds to the observation in \cref{ssec:1D} that, for higher order polynomials, the suppression of the highest eigenvalues eventually leads to one of the lower eigenvalues dominating.

\begin{figure}[!b]
    \centering
    \subfloat[$P=1$.]{\includegraphics[trim=80 70 80 70, clip, width=\linewidth]{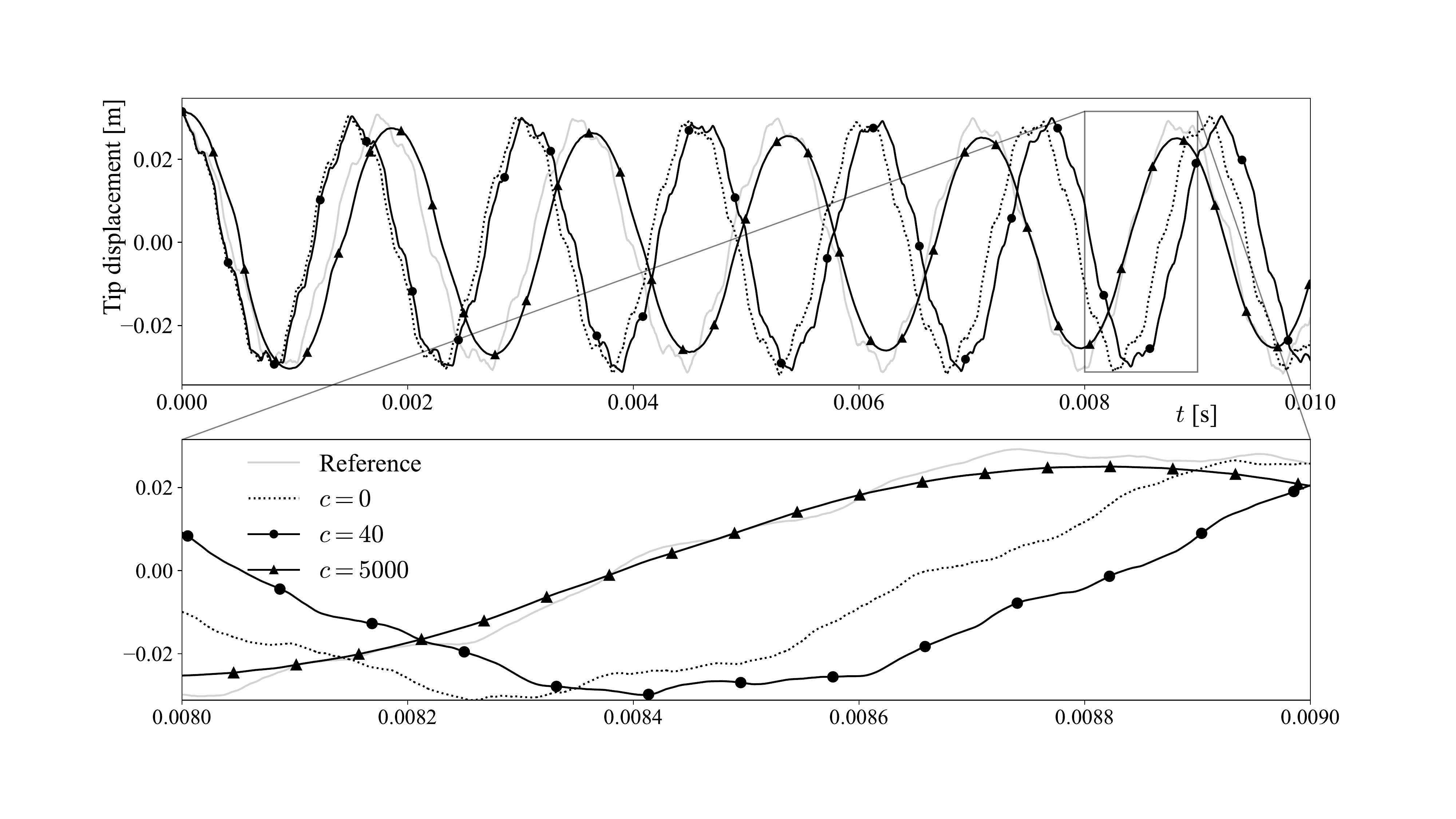} \label{fig:3Ddisp_P1} }\hspace{.25cm}
    \subfloat[$P=2$.]{\includegraphics[trim=80 70 80 50, clip, width=\linewidth]{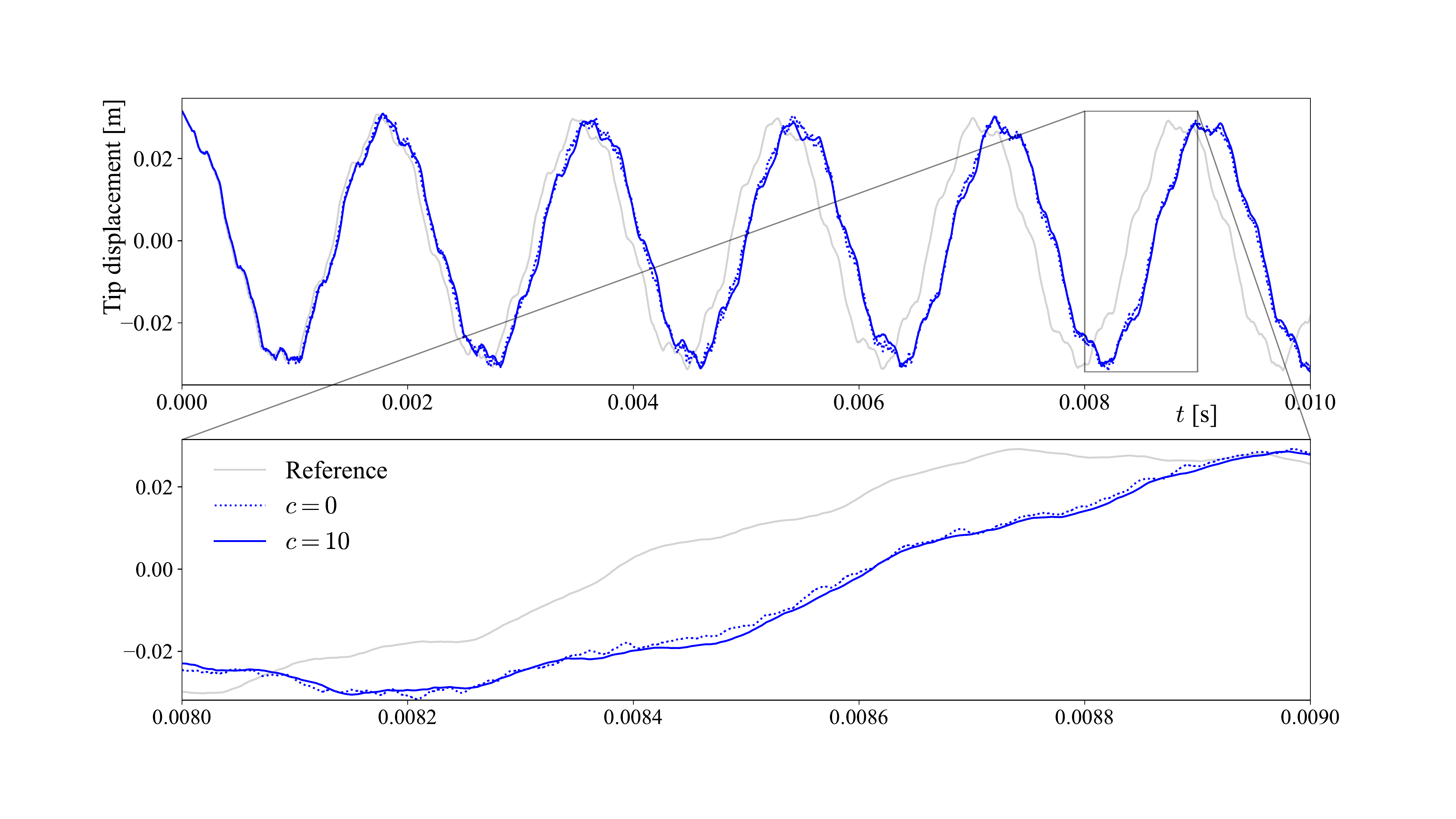} \label{fig:3Ddisp_P2} }
    \caption{Tip displacement.}\label{fig:3Ddisp}
\end{figure}

To study whether there is any quality loss when we make use of a $c$ that produces a significant cost-reduction, we choose more extreme values for $c$ than those in \cref{sec:dynamics}. For the computation with the linear basis functions we make use of $c=40$ and $c=5000$, which yield factors of increase of critical time step of approximately 6.5 and 22.2, respectively. For the computation with quadratic basis functions we use $c=10$, as this already roughly produces the plateau value of a factor 4.0 increase in critical time step. \Cref{fig:3Ddisp} shows the computed tip displacement as a function of time for the linear and quadratic cases computed on the coarse mesh of \cref{fig:3Dmesh1}. In both figures, the light-gray line is the reference solution computed with linear basis functions on the refined mesh that is depicted in \cref{fig:3Dmesh2}.

The results for both polynomial orders demonstrate that the impact of mass scaling with moderate choices of $c$ is minor compared to the source of error originating from the insufficient mesh resolution. In both cases we even observe small improvements. In \cref{fig:3Ddisp_P1}, the case of $c=0$ exhibits a severe overprediction of the dominant eigenfrequency. For $c=20$, this overprediction is reduced slightly while requiring 6.5 times less time step computations. The case $c=5000$ is chosen as an extreme as it shows that this behavior persists, and that the mass scaling can be tuned to yield the correct dominant frequency. Clearly, this comes at the cost of an apparent artificial damping, which reduces the accuracy of the high- and low-frequency response. %We anticipate that this artificial dampening is the result of the inconsistency of the new term, which adds resistance to rigid body rotation.

The quality of the coarse-mesh simulation is improved by the increase in polynomial order, as shown in \cref{fig:3Ddisp_P2}. For that simulation, the case of $c=5$ does not noticeably affect the dominant frequency within these first 5.5 periods. However, it does appear to suppress the artificial high-frequency oscillations that are visible in the dotted line in the zoom-in, and which are not present in the result of the reference computation.

\begin{figure}[!t]
    \centering
    \subfloat[$P=1$.]{\includegraphics[trim=12 10 20 0, clip, width=0.5\linewidth]{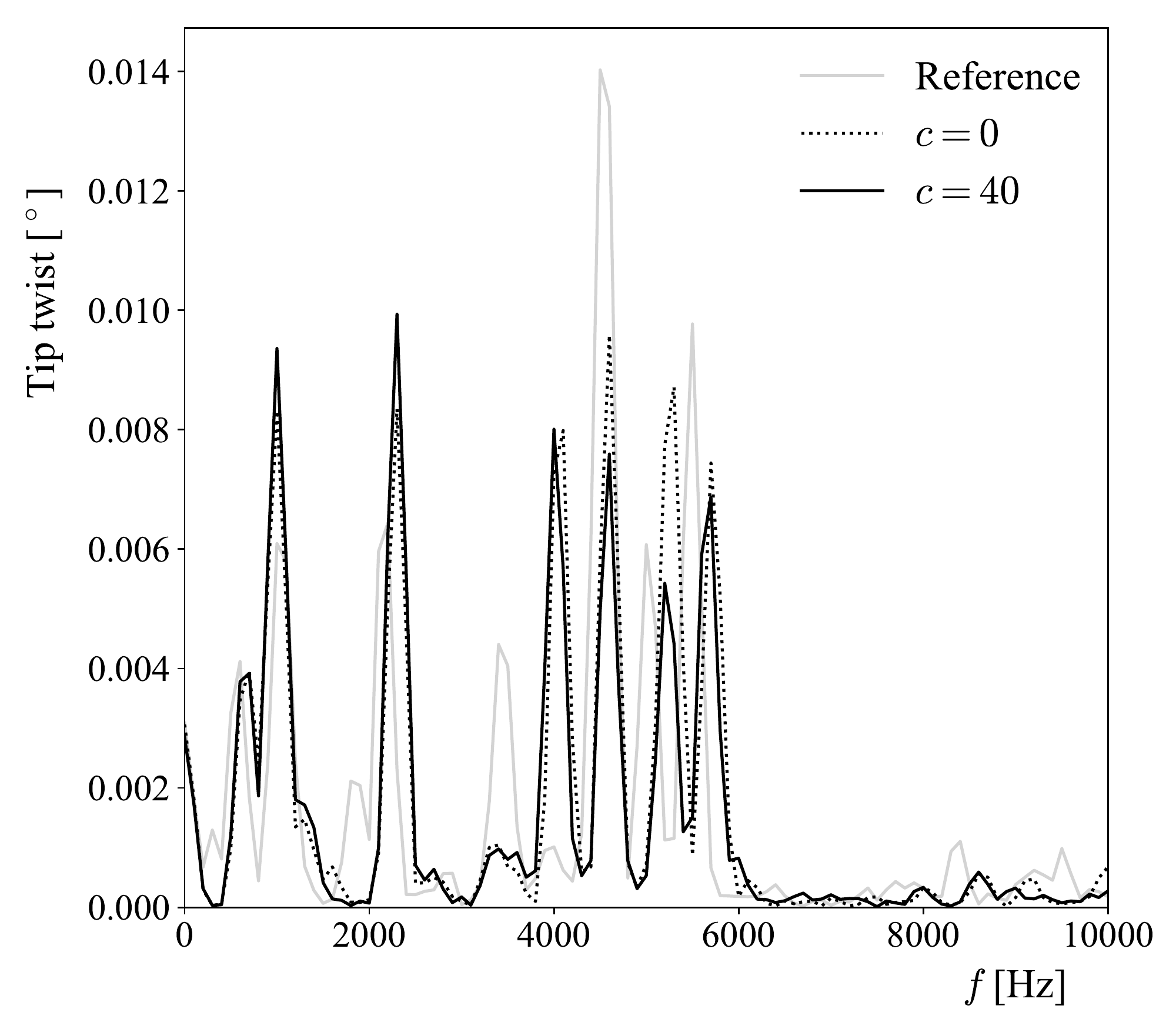} \label{fig:3Dfft_P1} }
    \subfloat[$P=2$.]{\includegraphics[trim=12 10 20 0, clip, width=0.5\linewidth]{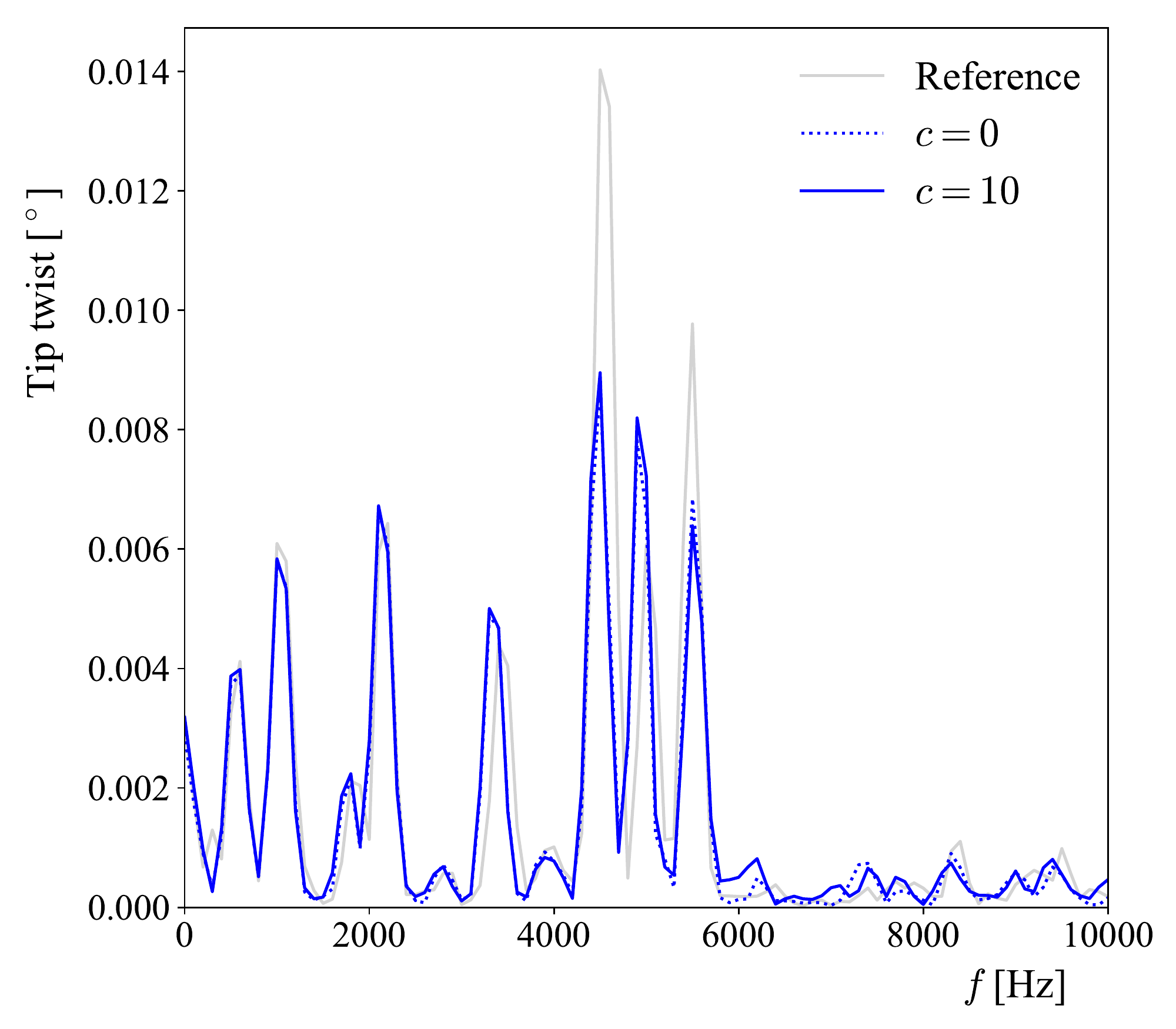} \label{fig:3Dfft_P2} }
    \caption{Fast Fourier transform of the angle of twist.}\label{fig:3Dfft}
\end{figure}

Due to the geometrical asymmetry, the beam response includes a component out of the $y$-plane. We compute the twist at the tip of the beam as a measure for this deformation. \Cref{fig:3Dfft} shows the frequency response thereof, as computed with a fast Fourier transform of the data from the first 0.01 seconds, windowed with a Blackman window \cite{Harris1978}. \Cref{fig:3Dfft_P1} shows the results for the computation with linear basis functions on the coarse mesh. The $c=5000$ case has been removed from the plot to avoid clutter and since it produces a largely inaccurate spectrum due to the artificial damping.
\Cref{fig:3Dfft_P2} shows the results for the computation with quadratic basis functions on the coarse mesh, which again results in a considerable improvement in the reproduction of the reference spectrum. For both computations, we observe that the scaled mass only marginally affects the system response. The locations and magnitudes of the peaks are not affected, and the changes in magnitudes are minor. For the linear basis functions, the overpredicted magnitude of one of the dominant peaks is even reduced. For the quadratic basis functions, the two spectra only become distinguishable past the 6000Hz frequency, at which range there are no more significant peaks.

\section{Conclusion and outlook}
\label{sec:conclusion}

In this article, we have proposed a variationally consistent addition to the mass matrix to suppress the highest eigenvalues in structural dynamics related applications. Our proposed term follows from a symmetric penalty of the natural transmission condition across element interfaces. For typical solid mechanics applications, the natural transmission condition is the interfacial equilibrium condition, which says that the jump of the traction should be zero. We weigh the new term by a small penalty, for which we derived an expression that ensures consistent results irrespective of the mesh size, polynomial order and spatial dimension.

For the linear wave equation, we have shown that the mass scaling improves the spectra on simple domains considerably, for all considered polynomial orders (up to quartics). In all cases, the low frequencies are virtually unaffected, the intermediate frequencies become more accurate, and the high frequencies are reduced. Also on more complicated domains and on irregular meshes, the mass-scaling term permits increases of critical time steps by factors of~2.5 without negatively affecting solution accuracy or convergence behavior.

We have also proposed a simplification that makes the mass-scaling term suitable for explicit time-stepping methods in non-linear solid mechanics applications. In our example problem, the mass scaling permits a factor 6.5 and factor 4.0 increase in the critical time-step size for linear and quadratic basis functions, respectively, without any negative impact on the solution accuracy (arguably, we even observe improvements).

There are various potential directions for future research. The proposed mass-scaling term may be extended to target the higher-order derivatives across element interfaces for higher-order polynomial basis functions, i.e., an application of the formulation from \cite{Nguyen2022a} for $C^0$-continuous basis functions. This would likely mitigate the plateauing behavior observed in \cref{ssec:1D,fig:3Dtcrit}. Secondly, one could investigate operator splitting techniques to move the mass-scaling term to the right-hand-side. This would guarantee the variational consistency of the formulation for a wide range of non-linear problems without the need for reassembly of the mass matrix in each time step. Additionally, operator splitting could provide a pathway for using our mass-scaling term in combination with a lumped mass matrix. %Lastly, our mass-scaling term could be added to the consistent mass matrix prior to lumping as a semi-consistent approach for obtaining a scaled lumped mass. Additionally, the penalty parameter can be chosen such that it outweighs the negative diagonal components due to row-summed mass lumping of higher-order serendipity elements.

\vspace{0.6cm}

\textbf{Acknowledgements:} The authors gratefully acknowledge support from the German Research Foundation (Deutsche Forschungsgemeinschaft) under grants SCH 1249/2-1 and SCH 1249/5-1.

\bibliographystyle{ieeetr}
%\bibliography{cm_references}
\bibliography{MyBib}

\end{document}